\documentclass{amsart}
\theoremstyle{plain}

\newcommand{\cleqn}{\setcounter{equation}{0}}
\newcommand{\clth}{\setcounter{theorem}{0}}
\newcommand {\sectionnew}[1]{\section{#1}\cleqn\clth}
\newcommand{\nn}{\hfill\nonumber}
\newtheorem{theorem}{Theorem}[section]
\newtheorem{lemma}[theorem]{Lemma}
\newtheorem{definition-theorem}[theorem]{Definition-Theorem}
\newtheorem{proposition}[theorem]{Proposition}
\newtheorem{corollary}[theorem]{Corollary}
\newtheorem{definition}[theorem]{Definition}
\newtheorem{example}[theorem]{Example}
\newtheorem{remark}[theorem]{Remark}
\newcommand \bth[1] { \begin{theorem}\label{t#1} }
\newcommand \ble[1] { \begin{lemma}\label{l#1} }

\newcommand \bpr[1] { \begin{proposition}\label{p#1} }
\newcommand \bco[1] { \begin{corollary}\label{c#1} }
\newcommand \bde[1] { \begin{definition}\label{d#1}\rm }
\newcommand \bex[1] { \begin{example}\label{e#1}\rm }
\newcommand \bre[1] { \begin{remark}\label{r#1}\rm }

\newcommand {\eth} { \end{theorem} }
\newcommand {\ele} { \end{lemma} }

\newcommand {\epr} { \end{proposition} }
\newcommand {\eco} { \end{corollary} }
\newcommand {\ede} { \end{definition} }
\newcommand {\eex} { \end{example} }  
\newcommand {\ere} { \end{remark} }

\newcommand \thref[1]{Theorem \ref{t#1}}
\newcommand \leref[1]{Lemma \ref{l#1}}
\newcommand \prref[1]{Proposition \ref{p#1}}

\newcommand \reref[1]{Remark \ref{r#1}}
\newcommand \lb[1]{\label{#1}}
\def \pab   {p^{\alpha,\beta}}
\def \pabn  {p^{\alpha,\beta}_n}
\def \Lab   {L_{\alpha,\beta}}
\def \Bab   {B_{\alpha,\beta}}
\def \pabe   {p^{\alpha,\beta}_\varepsilon}

\def \Labe   {L_{\alpha,\beta; \varepsilon}}
\def \Labet   {\widetilde{L}_{\alpha,\beta; \varepsilon}}
\def \pabet   {\tilde{p}^{\alpha,\beta}_\varepsilon}
\def \la {\lambda}
\def \lae {\lambda_\varepsilon}
\def \pz {\partial_z}
\def \px {\partial_x}
\def \a {\alpha}
\def \b {\beta}

\def \Ga {\Gamma}
\def \e {\varepsilon}
\def \ka {\kappa}
\def \Om {\Omega}
\def \Del {\Delta}
\def \th {\theta}
\def \ph {\varphi}
\def \Ph {\Phi}
\def \ps {\psi}
\def \Ps {\Psi}
\def \vph{\phi}
\def \Ii {{\mathcal I}}
\def \Con {{\mathcal C}}
\def \B {{\mathcal B}}
\def \K {{\mathcal K}}
\def \R {{\mathcal R}}
\def \D {{\mathcal D}^{(k, l)}_{\alpha, \beta; \varepsilon}}
\def \Dd {{\mathcal D}}
\def \Dko {{\mathcal D}^{(k, 0)}_{\alpha, \beta; \varepsilon}}
\def \Dol {{\mathcal D}^{(0, l)}_{\alpha, \beta; \varepsilon}}
\def \Bn {{\mathcal B}_{\alpha, \beta; \varepsilon}}
\def \Bz {{\mathcal B}'_{\alpha, \beta; \varepsilon}}
\def \Kn {{\mathcal K}_{\alpha, \beta; \varepsilon}}
\def \Kz {{\mathcal K}'_{\alpha, \beta; \varepsilon}}
\def \An {{\mathcal A}_{\alpha, \beta; \varepsilon}}
\def \Az {{\mathcal A}'_{\alpha, \beta; \varepsilon}}
\def \Rn {{\mathcal R}_{\alpha, \beta; \varepsilon}}
\def \Rz {{\mathcal R}'_{\alpha, \beta; \varepsilon}}
\def \Bnt {{\widetilde{\mathcal B}}_{\alpha, \beta; \varepsilon}}
\def \Bzt {{\widetilde{\mathcal B}}'_{\alpha, \beta; \varepsilon}}
\def \Knt {{\widetilde{\mathcal K}}_{\alpha, \beta; \varepsilon}}
\def \Kzt {{\widetilde{\mathcal K}}'_{\alpha, \beta; \varepsilon}}

\def \Rnt {{\widetilde{\mathcal R}}_{\alpha, \beta; \varepsilon}}
\def \Rzt {{\widetilde{\mathcal R}}'_{\alpha, \beta; \varepsilon}}
\def \Cset {{\mathbb C}}
\def \Zset {{\mathbb Z}}
\def \ra {\rightarrow}
\def \ol {\overline}
\def \wt {\widetilde}
\def \t {\tilde}
\def \go {\mapsto}
\def \sub { {\subset}}
\def \sup { {\supset}}

\def \det { {\mathrm{det}} }
\def \reg { {\mathrm{reg}} }
\def \Id { {\mathrm{Id}} }
\def \ad { {\mathrm{ad}} }
\def \Ad { {\mathrm{Ad}} }
\def \spa { {\mathrm{Span}} }
\def \ker { {\mathrm{Ker}} }
\begin{document}
\title[Discrete Bispectral Darboux transformations]
{Discrete Bispectral Darboux transformations from Jacobi operators}
\author[F.~A.~Gr\"unbaum]{F. Alberto Gr\"unbaum}\thanks{
The first author was partially
supported by NSF grant DMS94-00097}
\address{
Department of Mathematics \\
University of California at Berkeley \\
Berkeley, CA 94720, U.S.A.}
\email{grunbaum@math.berkeley.edu}
\author[M.~Yakimov]{Milen Yakimov}\thanks{
The second author was partially
supported by NSF grants
DMS94-00097 and DMS96-03239}
\address{
Department of Mathematics \\
University of California at Berkeley \\
Berkeley, CA 94720, U.S.A.}
\email{yakimov@math.berkeley.edu}
\date{October 2000}
\begin{abstract}
We construct families of bispectral difference operators of the 
form $a(n)T + b(n) + c(n) T^{-1}$ where $T$ is the shift operator.
They are obtained as discrete Darboux
transformations from appropriate extensions of Jacobi operators.
We conjecture that along with operators previously   
 constructed by Gr\"unbaum, Haine,
Horozov, and Iliev they exhaust all bispectral regular
(i.e. $a(n) \neq 0, c(n) \neq 0, \forall n \in \Zset)$ operators
of the form above.
\end{abstract}
\maketitle
\section{Introduction}
The bispectral problem, as originally  formulated by Duistermaat and
Gr\"unbaum
\cite{DG}, asks for a description of   
all situations where a pair of differential operators in the
variables
$x$ and $z$ have a common eigenfunction
$\Ps(x, z)$
\begin{align}
\lb{1}
&L(x, \px) \Ps(x, z) = \la(z) \Ps(x, z),
\\
\lb{2}
&B(z, \pz) \Ps(x, z) = \th(x) \Ps(x, z).  
\end{align}
For simplicity we say that $L$ , or $B$, or $\Ps$,  are bispectral
when the situation above holds.

The results in \cite{DG} already revealed a number of interesting
connections with a variety of topics ranging from the Korteweg--deVries
equation to the problem of isomonodromic deformations for differential
operators with rational coefficients.
Later even more unexpected
connections with different areas of pure mathematics were found.
These include automorphisms and ideal structure of the Weyl algebra
in one variable \cite{BW, BHY}, representations of the $W_{1+\infty}$
algebra \cite{BHYd}, Calogero--Moser system \cite{Wcm}, Huygens'
principle
\cite{B}, traces of intertwiners for representations of
(quantized) simple
Lie algebras 
\cite{EV, FMTV} (the last two in the multivariable case).

In \cite{DG} all bispectral differential
operators $L(x, \px)$ of second order were classified. Notice that
if one insists that $B(z, \pz)$ should also be of order two then one
is necessarily dealing with the Bessel or Airy cases.
In that paper very explicit use was made of the
Darboux transformation mapping a given second order
differential operator into another one. When starting from an
appropriate bispectral $L(x, \px)$ this was shown to produce another
such, with a different $B(z, \pz)$.
Wilson \cite{Wil}
approached the problem from the viewpoint of commutative algebras
of differential operators.
He classified
all maximal bispectral algebras of rank one (which by definition is the
greatest common divisor of the orders of all operators in the algebra).
In \cite{BHYc, KR} the idea of applying  Darboux transformations to
commutative algebras of differential operators was developed. This
allowed  for a unification of the apparently unrelated
methods in \cite{DG, Wil} and an extension of them to
the higher rank case. Further interesting results in this direction
were obtained in \cite{HM}.

Gr\"unbaum and Haine considered \cite{GH}
a discrete--differential version of the above problem
when the variable $x$ runs over the integer lattice $\Zset$
and accordingly one replaces
the differential operator $L(x, \px)$ by
a difference operator 
\[
L(n, T)= \sum_{i=p}^q b_i(n) T^i,
\quad b_p(n), b_q(n) \not\equiv 0
\]
acting on a function $f(n) \colon \Zset \ra \Cset$ by
\[
(L f)(n) = \sum_{i=p}^q b_i(n) f(n+i).
\]
Following \cite{vMM, Mum}, we define the support 
of $L(n,T)$ to be the ordered pair $[p,q].$ 
{\em{Such a difference operator will be called regular if
the first and the last coefficients $b_q(n)$ and $b_p(n)$
are nowhere vanishing functions on $\Zset.$}}

This problem is in a sense a
generalization of the problem of classifying
orthogonal polynomials which are eigenfunctions of differential
operators. The point is that the standard three term recursion
relation gives rise to a very special type of difference operator,
represented by a {\em{semiinfinite tridiagonal}} matrix.
In \cite{GH}, Gr\"unbaum and Haine showed that all instances
of difference operators with support $[-1,1]$ 
and second order differential operators satisfying
\eqref{1}--\eqref{2}
result by replacing the variable $n$ in the classical cases of
the Hermite, Laguerre, Jacobi, and
Bessel polynomials (discovered by Bochner, see \cite{B}) by a variable
$n+\e$ with n running over the integer lattice and $\e$ arbitrary
(see Sect.~2.2 for   
precise definitions). The differential operator in z is the celebrated
hypergeometric second order differential operator of Gauss.
It is worth noting that the corresponding
eigenfunctions $\Ps(n, z)$ are no longer polynomials.
For a recent survey of this area, see \cite{H2}.

Very recently Haine and Iliev classified all maximal bispectral
difference algebras of rank one \cite{HI}. Their result is a beautiful
extension
of Wilson's work \cite{Wil} where the Grassmannians associated to Darboux
transformations on differential operators are substituted with flag
varieties coming from such transformations on difference operators. Among
these some algebras that contain an operator with support
$[-1,1]$ were isolated in \cite{HI2}, where they were conjectured to be
all of this type.

The aim of this paper is to make progress in obtaining a 
discrete--continuous analog of the result of \cite{DG}, namely 
a classification of all discrete bispectral operators of the form
$a(n)T +b(n) + c(n) T^{-1}$ (referred to as the extended Bochner--Krall
problem in \cite{H2}). In \cite{GHH} the Darboux process was applied to
a biinfinite extension of the Laguerre difference operators
considered in \cite{GH}. A large class of bispectral difference 
operators of the form above was thus constructed and many properties 
of the resulting objects were analyzed in detail. It is fair to say 
that the results in \cite{GHH} provide a general treatment of the 
Laguerre case. The case of Jacobi difference operators has so far not 
been amenable to a similar  treatment and only some 
special cases of Darboux maps were proved to preserve the bispectral
property. The goal of the present paper is to provide such a 
general treatment in the Jacobi case and to state a conjecture
for the classification problem above. 

The rest of the introduction describes our results.

We take as a starting point the following natural extensions of the Jacobi
polynomials, constructed in
\cite{GH}:
\[
\pabe(n, z) = \frac{(\e+\a+1)_n}{(\e+1)_n}
              F(-(n+\e), n+\e+\a+\b+1, \a+1, (1-z)/2).
\]
Here and later we use $F$ for Gauss' 
${}_2 \hspace{-0.21mm} F_1$ hypergeometric function.
For negative integer values of $\a,$ see \eqref{1.16b}.
They are no longer polynomials but are still
eigenfunctions of a biinfinite difference
operator $\Labe(n,T)$ of the form $a_0(n) T + b_0(n) + c_0(n) T^{-1}$ and
a differential operator $\Bab(z, \pz):$
\begin{align*}
&\Labe(n, T) \pabe(n, z) = z \pabe(n, z),
\\
&\Bab(z, \pz) \pabe(n, z) = \lae(n) \pabe(n, z).
\end{align*}
The operator $\Labe(n, T)$ is obtained by the formal
change of variables $n \go n+\e$ from the standard (difference) 
Jacobi operator and is explicitly
defined in \eqref{1.17}. The operator $\Bab(z, \pz)$ and the 
spectral function $\lae(n)$ are given in
eqs. \eqref{1.8}, \eqref{lae}.

The sets of difference operators that we consider are obtained
by the following version of the Darboux map starting from the 
operators $\Labe(n, T).$ Let $P(n, T)$
be a regular difference operator whose kernel is preserved by 
$\Labe(n, T).$ Then there exists a (unique) difference
operator $L(n, T)$ such that
\begin{equation}
\lb{Dar}
L(n, T) P(n, T) = P(n, T) \Labe(n, T)
\end{equation}
which we refer to as {\em{a Darboux transformation from 
$\Labe(n,T).$}}
The advantage of this version is that this $L(n, T)$
is necessarily of the same form as $\Labe(n, T),$ i.e. 
$L(n, T) = a(n) T + b(n) + c(n) T^{-1}$ for some functions
$a(n),$ $b(n),$ $c(n),$ $n \in \Zset.$ 

If $q(x)$ denotes the characteristic polynomial of the endomorphism
$L(n, T)$ acting on the finite dimensional space $\ker P(n,T),$ the
\[
\ker P(n, T) \sub \ker q(\Labe(n,T)).
\]
In view of this it is natural to parametrize the sets
of operators $L(n,T)$ by the Grassmannians of special subspaces
of $\ker q(\Labe(n,T))$ that can occur as $\ker P(n, T).$
Denote the set of difference operators $L(n,T)$ corresponding
to characteristic polynomial $q(x)= (x-1)^k (x+1)^l$ by
\[
\D.
\] 
The operators in $\D$ are the main objects of study in this paper.
Their explicit form is given in Sect.~3.3.
Restricting to $q(x)$ with roots at $\pm 1$ guarantees that
$L(n,T)$ will have rational coefficients. This is an important 
feature of bispectral operators. See \cite{DG} in the 
differential case.

It is an easy consequence of \eqref{Dar} that
the function
\begin{equation}
\lb{Pss}
\Ps(n, z) = P(n, T) \pabe(n, z)
\end{equation}
is an eigenfunction of the operator $L(n, T),$ namely we have
\[  
L(n, T) \Ps(n, z) = z \Ps(n, z).
\]
Our main result is:
\bth{1.1} The difference operators $L(n,T)$ from the sets $\D$
are bispectral {\em{(}}or more precisely the functions 
$\Ps(n,z)$ \eqref{Pss} are eigenfunctions of differential 
operators in the variable $z)$ in the following cases

1.) $\a \in \Zset$ and $k \leq |\a|,$ $l=0,$

2.) $\b \in \Zset$ and $l \leq |\b|,$ $k=0,$

3.) $\a, \; \b \in \Zset$ and $k \leq |\a|,$ $l \leq |\b|.$
\eth

When $\e=0,$ and $k=1$ and (or) $l=1$ these results were established 
in \cite{KK, Zh}. All this work starts with the classical paper of
H.~L.~Krall, \cite{HLK}.

The proof of \thref{1.1} is based on a general result of
Bakalov, Horozov, and Yakimov \cite{BHY} which guarantees
that a Darboux transformation preserves the bispectral property
under some conditions on the operator $P(n, T).$
We will soon see that its application 
to the present situation is highly
nontrivial and requires, in particular, 
an intrinsic characterization of a space of
difference operators.

We will need some notation from \cite{BHY}, 
see Sect.~4.1 for more details. 
Denote by $\Bn$ the algebra of difference operators $S(n,T)$ with 
rational coefficients for which there exists a differential 
operator $G(z, \pz)$ (also having rational coefficients) 
satisfying
\begin{equation}
\lb{RSeq}
S(n, T) \pabe(n, z) = G(z, \pz) \pabe(n, z).
\end{equation}
All such operators $S(z, \pz)$ form a ``dual'' algebra $\Bz.$
The map 
\[
b \colon \Bn \ra \Bz, \quad
b( R(n, T) ) = S(z, \pz)
\]
is an antiisomorphism of associative algebras. 
Let $\Kn$ and $\Kz$ be the subalgebras of $\Bn$ and $\Bz$ consisting of
rational functions. Bispectrality of $\pabe(n, z)$ is equivalent
to $\Kn$ and $\Kz$ being both nontrivial. Finally we arrive at the most
important object for our consideration, namely the space
\begin{align*}
\Rn = \{ (\mu(n))^{-1} P_0(n, T) \, \mid & \; \mu(n) \in \Kn, \, 
P_0(n, T) \in \Bn, \text{and the operator} \\
&(\mu(n))^{-1} P_0(n, T) \, \text{does not have poles at} \, 
n \in \Zset \}.
\end{align*}
According to Theorem~1.2 of \cite{BHY}, $\Ps(n ,z)$ is an eigenfunction of
a differential operator in the variable $z,$ if 
\[
P(n, T) \in \Rn.
\]

Thus to prove \thref{1.1}
we need a {\em{good}} description of the space $\Kn$ which 
can be used to check whether the operators $P(n, T)$ from
\eqref{Dar} belong to $\Rn.$ This is the hardest step in our paper.
Let $\Del$ denote the algebra of abstract difference operators
with rational coefficients
of the form $\sum_{i=p}^{q} b_i(n) T^i$ with rational
functions $b_i(n)$ (possibly having poles in $\Zset).$
The key point of our approach is to consider the involution $I$ 
of $\Del$ acting on rational functions $h(n)$ by
\[
(I h) (n) := h(-(n+2 \e +\a +\b+1))
\]
and on the shift operator $T$ by $I(T) := T^{-1}.$
In Sect.~4.2 we prove that $\Rn$ consists of those 
difference operators from $\Del$ that do not have poles in $\Zset$ and
after conjugation with
the function
\[
\vph(n)= \frac{(\e+\a+1)_n}{(\e+1)_n}
\]
become $I$-invariant. 
 
The final step of the proof of \thref{1.1} is to show that
the hypothesis guarantee that the kernel of the operator 
$P(n,T)$ (defining $L(n,T))$ 
has a basis of functions $f(n)$ for which the ratio
$f(n)/\vph(n)$ is an (almost) $I$-invariant rational function in $n.$ 
This is done in Sect.~5.1. Finally Sect.~5.2 recapitulates the strategy of 
the proof of \thref{1.1} for the special case of the set 
$\Dd^{(2,0)}_{\a, \b, \e}.$ The reader may find it useful to consult
this section while reading the paper.

Let us also note that in the case of Laguerre
polynomials the situation simplifies a lot due to a presense 
of a relation of the type \eqref{RSeq} with a difference operator 
$S(n, T)$ of the form $s_1(n)T + s_0(n)$ and a first order differential
operator $G(z, \pz)$ (see expressions (2.3) and (2.8) in \cite{GHH}). 
It is not hard to show that as a consequence of this the
analog of $\Rn$ in that case is simply the space of difference operators
with rational coefficients.
 
Comparing with the differential case \cite{DG}, it is natural to
conjecture that {\em{all second order regular bispectral difference
operators {\em{(}}i.e. having support $[-1,1]${\em{\/)}} are
exhausted by the families of operators constructed in \cite{GH, GHH, HI2}
and in this article.}} The operators in 
\cite{HI2} are obtained as Darboux transformations from the operators
$\Labe(n, T)$ for half integer values of the parameters $\a, \b$ 
and are the analogs of ``KdV family'' in the differential case \cite{DG}.

For later use we introduce some convenient notation. If 
$f(n) \colon \Zset \ra \Cset$ is 
a nowhere vanishing function and $D_1(n, T),$ $D_2(n,T)$ are difference
operators we denote
\begin{align*}
& \Ad_{f(n)}D_1(n,T):= f(n) D_1(n, T) f(n)^{-1},
\\
& \ad_{D_2(n,T)}D_1(n,T):= D_2(n,T) D_1(n, T) - D_1(n, T) D_2(n,T).
\end{align*}   
\section{Biinfinite Jacobi operators}
In the first part of this section we review some properties of 
the classical Jacobi polynomials $\pabn(z).$ The second one discusses
certain functions $\pabe(n, z)$ which are eigenfunctions
of biinfinite analogs $\Labe(n,T)$ of the Jacobi difference operators. 
The third part describes Darboux maps between the operators $\Labe(n,T)$
with shifted indices $\a, \b.$ 
\subsection{Jacobi polynomials}
The Jacobi polynomials are the orthogonal polynomials for the
measure $(1-z)^\a (1+z)^\b dz$ on the interval $[-1, 1],$ 
$(\a, \b > -1),$ normalized by
\[
\pabn(1) = 2^{-n} \binom{n + \a}{n}, \; n \in \Zset_{\geq 0}.
\]
They are given by 
\begin{align} 
\lb{1.2}  
\pabn(z) 
&=\binom{n + \a}{n} F(-n, n+\a+\b+1; \a+1; 
                       (1-z)/2)
\end{align}
where $F(a,b;c;x)$ denotes the Gauss' hypergeometric function. 
The reader can consult \cite[pp. 209--217]{MOS}
for other explicit formulas and a list of major relations for 
$\pabn(z).$ 
Let
\begin{equation}
\lb{1.4}
\pab(n,z)=
  \begin{cases}
          \pab_n(z), &\text{for $n \in \Zset_{\geq 0}$} \\
          0,         &\text{for $n \in \Zset_{<0}$}
  \end{cases}
\end{equation}
Now $\pab(n, z)$ are functions of a discrete parameter $n$ and a
continuous parameter $z.$ They satisfy a three term recursion relation
\begin{equation}
\lb{1.6}
  \Lab(n, T) \pab(n, z) = z \pab(n, z)
\end{equation}
where $\Lab(n, T)$ are the difference operators
\begin{align}
\lb{1.7}
  \Lab(n, T)
  &=\frac{2(n+1)(n+\a+\b+1)}{(2n+\a+\b+1)(2n+\a+\b+2)} T
  +\frac{\b^2-\a^2}{(2n+\a+\b)(2n+\a+\b+2)}
\\
\nn
  &+\frac{2(n+\a)(n+\b)}{(2n+\a+\b)(2n+\a+\b+1)} T^{-1}
\end{align}
called {\em{Jacobi operators.}}
In addition $\pab(n, z)$ are eigenfunctions of the differential operators 
$\Bab(z, \pz)$ given by
\begin{equation}
\lb{1.8}
  \Bab(z, \pz) = (z^2-1) \pz^2+(\a-\b+(\a+\b+2)z)\pz,
\end{equation}
i.e.
\begin{equation}
\lb{1.9}
  \Bab(z, \pz) \pab(n, z) =\la(n) \pab(n, z)
\end{equation}
for
\begin{equation}
\lb{1.10}
  \la(n) =n(n +\a + \b+ 1).
\end{equation}
In view of \eqref{1.6} and \eqref{1.9}, $\pab(n, z)$ are 
discrete--continuous bispectral functions and $\Lab(n,T),$ $B(z, \pz)$
bispectral difference (differential) operators.

\subsection{The functions $\pabe(n,z)$}
In a study of the relation between the so called ``associated Jacobi
polynomials'' and the discrete--continuous bispectral problem 
Gr\"unbaum and Haine, see \cite{GH, GH2, H}, introduced
the functions
\begin{equation}
\lb{1.16}
\pabe(n, z) = \frac{(\e+\a+1)_n}{(\e+1)_n}
              F(-(n+\e), n+\e+\a+\b+1; \a+1; (1-z)/2)
\end{equation}
$(n \in \Zset, z \in \Cset, |z|<1)$ defined for those 
$\e, \a, \b \in \Cset$ such that $\a \notin \Zset_{<0},$ 
and $\e \notin \Zset_{<0},$ $\e+\a \notin \Zset_{\geq 0}.$
We will see later that the first restriction can be lifted.

The functions $\pabe(n, z)$ are no longer polynomials
but satisfy relations, similar to the ones for $\pab(n, z).$
In particular they are eigenfunctions of the following
difference operators with support $[-1, 1]$
\begin{align}
\lb{1.17}
\Labe(n, T)&=\frac{2(n+\e+1)(n+\e+\a+\b+1)}
           {(2n+2\e+\a+\b+1)(2n+2\e+\a+\b+2)} T
\\
\nn
     &+\frac{\b^2-\a^2}
        {(2n+2\e+\a+\b)(2n+2\e+\a+\b+2)}+
\\
\nn 
     &+\frac{2(n+\e+\a)(n+\e+\b)}
        {(2n+2\e+\a+\b)(2n+2\e+\a+\b+1)} T^{-1}
\end{align}
{\em{and}} of the differential operators $\Bab(z, \pz),$ eq. \eqref{1.8}. 
The corresponding relations are
\begin{align}
\lb{1.19} 
&\Labe(n, T) \pabe(n, z) = z \pabe(n, z),
\\
\lb{1.20}
&\Bab(z, \pz) \pabe(n, z) = \lae(n) \pabe(n, z).
\end{align}
where
\begin{equation}
\lb{lae}
\lae(n) = (n+\e)(n + \e + \a + \b + 1).
\end{equation}

The difference operators $\Labe(n, T)$ will still be called
{\em{Jacobi operators.}} Further we will only deal with the case when they
are regular, i.e. when their coefficients of $T$ and $T^{-1}$ do not
vanish for $n \in \Zset.$ This amounts  to the conditions
\begin{equation}
\lb{cond}
\e, \e+\a, \e+ \b, \e + \a + \b, 2 \e+ \a +\b \notin \Zset.
\end{equation}
{\em{It may be useful to stress here that these will eventually be the 
only  restrictions on our parameters $\a,$ $\b,$ $\e.$}}

The operators $\Labe(n, T)$ do satisfy certain ``transformation 
properties''. For instance the following relations hold
\begin{align}
\lb{rel1}
&L_{-\a, -\b, \e+\a+\b}(n, T) = \Labe(n, T), \\
\lb{rel2}
&\Ad_{(-1)^n} L_{\b, \a, \e}(n, T) =
 \Ad_{(-1)^n} L_{-\b, -\a, \e+\a+\b}(n, T)= - \Labe(n, T).
\end{align} 
It is tempting to use \eqref{rel1} to limit attention to the case
$\a \geq 0.$ However, this would eventually bring an undesirable degree
of asymmetry in the treatment of the parameters $\a$ and $\b.$
For this reason we prefer to introduce the appropriate functions
$\pabe(n,z)$ for $\a \in \Zset_{<0}$
(and $\e, \e+\b, \e+\a, \e+\a+\b \notin \Zset)$ by using
\eqref{1.16} and recalling, see \cite[p. 38]{MOS}
that for $m \in \Zset_{\geq 0}$
\begin{equation}
\lb{*}
\lim_{c \ra -m} 
\frac{1}{\Ga(c)} F(a, b;c;z)=
\frac{(a)_{m+1} (b)_{m+1}}{(m+1)!} z^{m+1} 
F(a+m+1, b+m+1; m+2;z).
\end{equation}
We see below that this leads to the following expression for 
$\pabe(n, z)$ with $\a \in \Zset_{<0}$ (as long as 
\eqref{cond} is satisfied) 
\begin{equation}
\lb{1.16b}
    \Con
    \frac{(\e+\b+1)_n}{(\e+\a+\b+1)_n}
    \frac{(1-z)^{-\a}}{2^{-\a}}  
  F(-(n+\e+\a), n+\e+\b+1;-\a+1;(1-z)/2)
\end{equation}
where the constant $\Con=\Con(\a, \b, \e)$ is explicitly given by
\[
\Con=\Con(\a, \b, \e)= \frac{(-1)^\a}{(-\a-1)!} \cdot
\frac{(-\e)_{-\a}(\e+\a+\b+1)_{-\a}}{(-\a)!}.
\]
It is easy to check that the assumptions \eqref{cond} imply that
$\Con(\a, \b, \e)$ is well defined and does not vanish.

The expression above can be derived by a continuity argument
using \eqref{1.16} and \eqref{*}
when $\a$ approaches a value in $\Zset_{<0}.$ 
To see this it is important to notice that
for $\a \in \Zset_{<0}$ the identities 
\[
\frac{(\e+\a+1)_n}{(\e+1)_n}=
\frac{(-\e)_{-\a}}{(-(n+\e))_{-\a}} \quad 
\text{and} \quad
\frac{(\e+\b+1)_n}{(\e+\a+\b+1)_n}=
\frac{(n+\e+\a+\b+1)_{-\a}}{(\e+\a+\b+1)_{-\a}}
\]
allow one to rewrite the factor
\[
\frac{(\e+\a+1)_n}{(\e+1)_n} \cdot
\frac{(-(n+\e))_{-\a} (n+\e+\a+\b+1)_{-\a}}{(-\a)!}
\]
as
\begin{equation}
\lb{factor}
\frac{(-\a-1)!}{(-1)^\a} \cdot 
\Con(\a, \b, \e) \cdot
\frac{(\e+\b+1)_n}{(\e+\a+\b+1)_n}
\end{equation}
which except for the first constant is the factor in front 
of \eqref{1.16b}.

Then conditions \eqref{cond} guarantee that $\pabe(n,z)$ is well defined
(see \eqref{1.16} and \eqref{1.16b}) and satisfies \eqref{1.19}
and \eqref{1.20}. It was proved in \cite{GH} that the space of common
solutions of \eqref{1.19} and \eqref{1.20} in a domain $\Om \sub \Cset,$   
not containing $\pm 1,$ is two dimensional.
Notice also that \eqref{cond} excludes, in particular, the original
operators $\Lab(n, z)$ $(\e =0)$ since their leading coefficient
vanishes for $n=-1.$

Finally we explain how \eqref{1.19} follows from \eqref{1.6}.  
Conjugate the operator $\Lab(n,T)$ with 
$\binom{n+a}{n}=\frac{(\a+1)_n}{n!}.$ The resulting difference
operator has rational coefficients and the eigenfunction 
$F(-n, n+\a+\b+1; \a+1; (1-z)/2),$ cf. \eqref{1.2}. 
The operator obtained from it by the formal change $n \go n+\e$ has
the eigenfunction $F(-n-\e, n + \e + \a+\b+1; \a+1; (1-z)/2)$
and all we need to do is conjugate it with $(\e+1)_n/(\e+\a+1)_n.$ The
result is the operator $\Labe(n,T)$ which proves \eqref{1.19} 
in the case $\a \notin \Zset_{<0}.$ The case $\a \in \Zset_{<0}$
follows from the definition \eqref{1.16} using the limit
\eqref{*}. 
\subsection{Darboux maps between Jacobi operators}
There are four difference relations connecting the values of the Jacobi
polynomials $\pab(n, z)$ with shifted indices:
\begin{align*}
& p^{\a-1, \b}(n,z)= \left(
       \frac{n+\a+\b}{2n +\a +\b} -   
       \frac{n+\b}{2n +\a +\b} T^{-1} \right) \pab(n,z),
\\
& \pab (n,z)= \frac{1}{z-1} \left(
       \frac{2(n+1)}{2n +\a +\b+1} T -
       \frac{2(n+\a)}{2n +\a +\b+1} \right) p^{\a-1, \b}(n,z),
\end{align*}
and
\begin{align*}
& p^{\a, \b-1}(n,z)= \left(
       \frac{n+\a+\b}{2n +\a +\b} +
       \frac{n+\a}{2n +\a +\b} T^{-1} \right) \pab(n,z),
\\
& \pab(n,z)= \frac{1}{z+1} \left(
       \frac{2(n+1)}{2n +\a +\b+1} T +
       \frac{2(n+\b)}{2n +\a +\b+1} \right) p^{\a, \b-1}(n,z),
\end{align*}
(see for instance, \cite[eqs. pp. 209--219]{MOS}). 
Similarly to the proof of \eqref{1.19} at the end of the previous
subsection, one shows the
following analogs of these identities for $\pabe(n, z)$
\begin{align}
\lb{1.22} 
&p_\e^{\a-1, \b}(n, z) = D_{-}^\a(n, T) \pabe(n, z), \; 
p_\e^{\a+1, \b}(n, z) = \frac{1}{z-1} D_{+}^\a(n, T) \pabe(n, z),
\\
\lb{1.24} 
&p_\e^{\a, \b-1}(n, z) = D_{-}^\b(n, T) \pabe(n, z), \; 
p_\e^{\a, \b+1}(n, z) = \frac{1}{z+1}D_{+}^\b(n, T) \pabe(n, z),
\end{align}
where the operators $D_\pm^\a(n, T)$ and $D_\pm^\b(n, T)$ 
are given by
\begin{align*}
& D_{-}^\a(n, T)= \left( \frac{\e+\a}{\a} \right)
  \left( \frac{n+\e+\a+\b}{2n+2\e+\a+\b} -
  \frac{n+\e+\b}{2n+2\e+\a+\b} T^{-1} \right),
\\
& D_{+}^\a(n, T)= \left( \frac{\a+1}{\e+\a+1} \right)
  \left( \frac{2(n+\e+1)}{2n+2\e+\a+\b+2} T -
  \frac{2(n+\e+\a+1)}{2n+2\e+\a+\b+2} \right),
\\
& D_{-}^\b(n, T)=
  \left( \frac{n+\e+\a+\b}{2n+2\e+\a+\b} +
  \frac{n+\e+\a}{2n+2\e+\a+\b} T^{-1} \right),   
\\
& D_{+}^\b(n, T)=
  \left( \frac{2(n+\e+1)}{2n+2\e+\a+\b+2} T +
  \frac{2(n+\e+\b+1)}{2n+2\e+\a+\b+2} \right).      
\end{align*}
The constant in \eqref{1.16b} was chosen to make the relations
\eqref{1.22}--\eqref{1.24} hold for all $\a \in \Cset.$
We show only the dependence on the index $\a$ of the operators
$D_{\pm}^\a(n, T)$ because the index $\b$ is unchanged in both sides of
eqs. \eqref{1.22}, similarly for the operators 
$D_{\pm}^\b(n, T).$ Eqs. \eqref{1.22}--\eqref{1.24} and
\eqref{1.19} imply the following factorizations
\begin{align}
& \Labe(n, T) - 1= 
D_{+}^{\a-1}(n, T) D_{-}^\a(n, T) = 
D_{-}^{\a+1}(n, T) D_{+}^\a(n, T),
\\
& \Labe(n, T) + 1= 
D_{+}^{\b-1}(n, T) D_{-}^\b(n, T) = 
D_{-}^{\b+1}(n, T) D_{+}^\b(n, T).
\end{align}
Hence the operators $L_{\a \pm 1, \b;\e}(n, T),$ 
$L_{\a, \b \pm 1 ;\e}(n, T)$ are Darboux transformations from
$\Labe(n, T)$ and eqs. \eqref{1.22}, \eqref{1.24}
represent the Darboux maps $\pabe(n,z) \go p^{\a \mp 1, \b}_\e(n,z)$
and $\pabe(n, z) \go p^{\a, \b \mp 1}_\e(n,z).$
\sectionnew{Darboux transformations from Jacobi operators}
The first part of this section contains some general facts about
discrete Darboux transformations in the form in which they will be 
used later (see, for instance, \cite{W} for the differential case).
The goal of the second part is an explicit description of
the kernels of the operators $(\Labe-1)^k (\Labe+1)^l.$ 
Based on it, in the third part we construct Darboux
transformations from $\Labe(n,T)$ which are the main objects of 
study in the rest of the paper. The conditions \eqref{cond} are 
assumed throughout Sects.~3.2--3.3.
\subsection{General remarks on Darboux transformations}
One says that the difference operator $L(n, T)$ is obtained by a Darboux
transformation from the difference operator $L_0(n, T)$ if there exists an
operator $P(n, T)$ such that
\begin{equation}
\lb{darboux}
L(n, T) P(n, T) = P(n, T) L_0(n, T).
\end{equation}
Assume that $L_0(n, T)$ has an eigenfunction $\Ps_0(n, z),$ i.e.
\begin{equation}
\lb{0eig}
L_0(n, T) \Ps_0(n, T) = g_0(z) \Ps_0(n)  
\end{equation}
for some function $g_0(z).$ Then
\[
\Ps(n, z):= P(n, T) \Ps_0(n, z)
\]
is an eigenfunction of $L(n, T):$
\begin{equation}
\lb{leig}
L(n, T) \Ps(n, T) = g_0(z) \Ps(n).
\end{equation}
The map $\Ps_0(n, T) \go \Ps(n, T)$ is also called a Darboux 
transformation. 

{\em{An important feature of the transformation \eqref{darboux}
for a regular difference operator $P(n, T)$ is that the operator
$L(n, T)$ has the same support as $L_0(n, T).$ Besides this
$L(n, T)$ is regular if and only if $L_0(n,T)$ is regular.}}

Given a difference operator $L_0(n, T),$ all
transformations of the type \eqref{darboux} 
with a regular difference operator $P(n, T)$ 
can be described in terms of the kernel of $P(n, T).$

\bpr{darb} {\em{(i)}} For a regular difference operator $P(n, T)$ there
exists
an operator $L(n, T)$ for which \eqref{darboux} holds if and only if
\begin{equation}
\lb{ker-darboux}
L_0(n,T) ( \ker P(n, T) ) \sub \ker P(n, T).
\end{equation}
The operator $L(n,T)$ satisfying \eqref{darboux} is unique.
    
{\em{(ii)}} Let $P(n, T)$ be a regular difference operator
satisfying \eqref{ker-darboux} and $q(x)$ be the characteristic 
polynomial of the linear map $L_0(n, T)$ acting in the 
space $\ker P(n, T).$ Then $\ker P(n, T) \sub q(L_0(n, T))$ 
and there exists an operator $Q(n, T)$ such that
\begin{align}
\lb{QP}
&q( L_0(n, T) ) = Q(n, T) P(n, T),
\\
\lb{PQ} 
&q(L(n, T) ) =  P(n, T) Q(n, T).
\end{align}
\epr

Note that the kernel of a regular difference operator $P(n, T)$ 
is finite dimensional. More precisely, if $P(n, T)$ has support
$[m_1, m_2]$ for some $m_i \in \Zset,$ 
then $\dim \ker P(n, T) = m_2 - m_1.$ For any $j \in \Zset$
the map 
\begin{equation}
\lb{ker-som}
f \go (f(j+1), \ldots, f(j+m_2 - m_1)), \; \mbox{for} \; 
   f: \Zset \ra \Cset 
\end{equation}
provides an isomorphism between $\ker P(n, T)$ and $\Cset^{m_2 -m_1}.$ 

The transformation $Q(n, T) P(n, T) \go P(n, T) Q(n, T)$ is a more
traditional version of the Darboux map. Although it is a special case
of the transformation $L_0(n, T) \go L(n, T)$ from eq. \eqref{darboux}, 
\prref{darb} shows that there always exists a polynomial $q(x)$
for which $q(L_0(n, T)) \go q(L(n, T))$ is a Darboux map in this sense.
\\
\noindent
{\em{Proof of \prref{darb}.}} (i) If $P(n, T),$ $L(n, T)$ satisfy
\eqref{darboux} and $f(n) \in \ker P(n, T)$ then
\[
P(n, T) (L_0(n, T) f(n) ) =
L(n, T) P(n, T) f(n) =0
\]
which proves \eqref{ker-darboux}. 

In the opposite direction, let us notice that a comparison of the
coefficients of the two sides of eq. \eqref{darboux} for a fixed
value of $n$ gives a finite system for the corresponding coefficients of
the unknown operator $L(n, T)$ having the same support as $L_0(n, T).$ 
One shows that it has a solution using
the standard linear algebra fact that for a finite matrix $A$ 
the system $A u = b$ has a solution if
and only if $v^t b = 0,$ $\forall v \in \ker A^t.$
In the particular case which we consider the last condition is fulfilled
because
of \eqref{ker-darboux}. 

The regularity of the difference operator $P(n,T)$ implies the uniqueness
of the operator $L(n, T)$ satisfying \eqref{darboux}. Indeed if there are
two such operators $L(n,T)$ and $L'(n, T)$ one can subtract the resulting
equalities \eqref{darboux}. This gives $(L(n,T) - L'(n, T)) P(n, T)=0$ 
which is a contradiction.

(ii) The relation $\ker P(n, T) \sub q(L_0(n, T))$ follows from the
definition of $q(x).$ Similarly to part (i), this implies the existence of
an operator $Q(n, T)$ satisfying \eqref{QP}. Eqs. \eqref{darboux} and
\eqref{QP} imply
\[
q(L(n, T)) P(n, T) = P(n, T) q(L_0(n, T)) =
(P(n, T) Q(n, T) ) P(n, T)
\]
and as a consequence of this \eqref{PQ}.
\hfill \qed 

A regular difference operator is reconstructed from its kernel
by the following lemma.

\ble{regdiff} Assume that $P(n, T)$ is a regular difference operator
with support $[m_1, m_2]$ and leading coefficient 1.
Let $\ker P(n, T) = \spa \{f^{(i)}(n)\}_{i=1}^m$ where $m=m_2 -m_2.$
Then the function
\[
\det(n):= \det( f^{(i)}(n-j) )_{i, j=1, m_1}^{m, m_2-1}
\]
does not vanish for $n \in \Zset$  and
\begin{equation}
\lb{operP}
P(n, T) = \frac{1}{\det(n)}
\begin{vmatrix}
  f^{(1)}(n+m_1) & \cdots &
  f^{(m)}(n+m_1) &
  T^{m_1} \cr
\cdots & \cdots & \cdots &
  \cdots \cr
  f^{(1)}(n+m_2) & \cdots &
  f^{(m)}(n+m_2) & 
  T^{m_2}
\end{vmatrix}
\end{equation}
where the determinant is expanded from left to right {\em{(}}the shift
operator $T$ does not commute with function multiplication{\em{\/).}}
\ele
\begin{proof}
The fact that the map \eqref{ker-som} is an isomorphism between 
$\ker P(n, T)$ and $\Cset^m$ implies that $\det(n)$ does not
vanish for $n \in \Zset.$
Clearly the functions $f^{(i)}(n)$ belong to the kernel 
of the operator in the r.h.s. of \eqref{operP}. It has
leading term 1 and the nonvanishing of $\det(n)$ implies \eqref{operP}. 
\end{proof}
\bre{comp}The composition of two Darboux transformations 
$L_0(n, T) \go L_1(n, T)$ and $L_1(n, T) \go L_2(n, T)$
of the type \eqref{darboux} is Darboux transformation
$L_0(n,T) \go L_2(n,T)$ of the same type. Indeed if
\[
L_i(n,T) P_i(n,T) = P_i(n,T) L_{i-1}(n,T), \,
i=1,2,
\]
then
\[
L_2(n,T) P_2(n,T) P_1(n,T)=
P_2(n,T) P_1(n,T) L_0(n,T).
\]
\ere
\subsection{Description of $\ker (\Labe-1)^k(\Labe+1)^l$} 
The main idea is to first find some functions $\ph(n, z)$ 
(depending on $\a, \b,$ and $\e$) such that
\begin{equation}
\lb{2.11}
\Labe(n, T) \ph(n, z) = z \ph(n, z)
\end{equation}
and then to consider the derivatives
\[
\ph_{\pm}^{(i)}(n) = \frac{1}{i!} \pz^i \ph(n, z) \Big|_{z = \pm 1}, 
\; i \in \Zset_{\geq 0}.
\] 
They satisfy
\begin{equation}
\lb{2.12}
(\Labe(n, T) \mp 1) \ph_{\pm}^{(i)}(n) = \ph_{\pm}^{(i-1)}(n), \;
\forall \; i \in \Zset_{\geq 0}
\end{equation}
with $\ph_{\pm}^{(-1)}(n) = 0.$ As a consequence of this
\[
(\Labe(n, T) \mp 1)^i \ph_{\pm}^{(j)}(n) =0, \;
\forall \; i \in \Zset_{> 0}, \; j=0, \ldots, i-1.  
\]

Before stating the results from this subsection we recall a relation
for the hypergeometric function that is a consequence of Gauss' 
relations
between contiguous hypergeometric functions.  Denote 
$F = F(a, b; c;(1-z)/2),$ $T F = F(a-1, b+1; c;(1-z)/2),$ 
and $T^{-1}F = F(a+1, b-1; c;(1-z)/2).$ Then for 
$c \notin \Zset_{\leq 0}$
\begin{multline}
\lb{2.14}
\frac{2(c-a)b}{(b-a)(b-a+1)} T F +
        \frac{2(a+b-1)(-2c+a+b+1)}{(b-a-1)(b-a+1)} F 
\\
+ \frac{2a(c-b)}{(b-a)(b-a-1)} T^{-1} F = z F.
\end{multline}
This can also be checked directly using the standard expansion
of $F(a, b;c, x)$ for $|x| <1,$ $c \notin \Zset_{\leq 0}$ 
\begin{equation}
\lb{expansion}
F(a, b; c; x) = \sum_{j=0}^\infty
\frac{(a)_j(b)_j}{j! (c)_j} x^j.
\end{equation}

\ble{2.4} The four functions
\begin{align}
\lb{2.15} 
&   \ph_{+}(n, z) =
    \frac{(\e+\a+1)_n}{(\e+1)_n}
    F(-(n+\e), n+\e+\a+\b+1;\a+1;(1-z)/2),
\\
\lb{2.16}
&   \ps_{+}(n, z) =
    \frac{(\e+\b+1)_n}{(\e+\a+\b+1)_n}
    F(-(n+\e+\a+\b), n+\e+1;-\a+1;(1-z)/2),
\\
\lb{2.17}
&   \ph_{-}(n, z) = 
    \frac{(-1)^n (\e+\b+1)_n}{(\e+1)_n}
    F(-(n+\e), n+\e+\a+\b+1;\b+1;(1+z)/2),
\\
\lb{2.18}   
&   \ps_{-}(n, z) = 
    \frac{(-1)^n(\e+\a+1)_n}{(\e+\a+\b+1)_n}
    F(-(n+\e+\a+\b), n+\e+1;-\b+1;(1+z)/2)
\end{align}
satisfy
\begin{equation}
\lb{2.19} 
(\Labe(n, T) - z) \ph_{\pm}(n, z) =
(\Labe(n, T) - z) \ps_{\pm}(n, z) =0,
\end{equation}
provided that $\a \notin \Zset_{<0} \, (\Zset_{>0})$ for
$\ph_+(n, z)$ $(\ps_+(n, z))$ and 
$\b \notin \Zset_{<0} \, (\Zset_{>0})$ for
$\ph_-(n, z)$ $(\ps_-(n, z)).$
\ele

Note that the assumptions \eqref{cond} guarantee that the denominators of
the
first factors of $\ph_\pm(n, z)$ and $\ps_\pm(n, z)$ do not vanish. 
\begin{proof}
The relation \eqref{2.19} for $\ph_+(n, z)$ holds
because $\ph_{+}(n, z) = \pabe(n, z).$ 
To check the one for $\ps_+(n, z),$ we conjugate 
$\Labe(n, T)$ by $(\e+\b+1)_n/(\e+\a+\b+1)_n$ 
(the factor in front of the r.h.s. of 
\eqref{2.16}).

The result
is 
\begin{align*}
& \Ad_{(\e+\b+1)_n/(\e+\a+\b+1)_n} \Labe(n, T)
\\
& =\frac{2(n+\e+1)(n+\e+\b+1)}
  {(2n+2\e+\a+\b+1)(2n+2\e+\a+\b+2)} T 
\\
& +\frac{\b^2-\a^2}
  {(2n+2\e+\a+\b)(2n+2\e+\a+\b+2)} 
\\
& +\frac{2(n+\e+\a)(n+\e+\a+\b)}
  {(2n+2\e+\a+\b)(2n+2\e+\a+\b+1)} T^{-1}. 
\end{align*}
This is the difference operator from the l.h.s. of \eqref{2.14}
with $a =-(n+\e+\a+\b),$ $b=n+\e+1,$ and $c =-\a+1$ which gives the
proof of \eqref{2.19} for $\ps_+(n,z).$  
The cases of $\ph_-(n, z)$ and $\ps_-(n, z)$ are 
handled in a similar fashion. 
\end{proof}

Next we consider the derivatives of 
$\ph_{+}(n, z),$ $\ps_{+}(n, z)$ at $z=1$ and of 
$\ph_{-}(n, z),$ $\ps_{-}(n, z)$ at $z=-1:$
\begin{align*}
& \ph_{\pm}^{(i)}(n) := \frac{1}{i!} \pz^i \ph_{\pm}(n, z)
  \Big|_{z=\pm1},
\\
& \ps_{\pm}^{(i)}(n) := \frac{1}{i!} \pz^i \ps_{\pm}(n, z)
  \Big|_{z=\pm1},
\end{align*}
$i \in \Zset_{\geq 0}$ (with the restrictions on $\a$ and $\b$ made at 
the end of \leref{2.4}).

Using the expansion \eqref{expansion} of the hypergeometric function, 
we obtain the following explicit formulas for $\ph_{\pm}^{(i)}(n)$
and $\ps_{\pm}^{(i)}(n)$
\begin{align}
\lb{2.25} 
&  \ph_{+}^{(i)}(n) = \frac{(\e+\a+1)_n}{(\e+1)_n} \cdot
   \frac{(-(n+\e))_i(n+\e+\a+\b+1)_i}{(-2)^i i! (\a+1)_i}
\\
\lb{2.26} 
&  \ps_{+}^{(i)}(n) = \frac{(\e+\b+1)_n}{(\e+\a+\b+1)_n} \cdot
   \frac{(-(n+\e+\a+\b))_i(n+\e+1)_i}{(-2)^i i! (-\a+1)_i}
\\
\lb{2.27} 
&  \ph_{-}^{(i)}(n) = \frac{(-1)^n(\e+\b+1)_n}{(\e+1)_n} \cdot
   \frac{(-(n+\e))_i(n+\e+\a+\b+1)_i}{2^i i! (\b+1)_i}
\\
\lb{2.28}   
&  \ps_{-}^{(i)}(n) = \frac{(-1)^n (\e+\a+1)_n}{(\e+\a+\b+1)_n} \cdot
   \frac{(-(n+\e+\a+\b))_i (n+\e+1)_i}{2^i i! (-\b+1)_i} \cdot
\end{align}
We define $\ph_{+}^{(i)}(n)$ $(\ps_{+}^{(i)}(n))$ for
$\a \in \Zset_{<0}$ $(\a \in \Zset_{>0}),$ $i<|\a|$ by 
\eqref{2.25}, \eqref{2.26} and 
$\ph_{-}^{(i)}(n)$ $(\ps_{-}^{(i)}(n))$ for
$\b \in \Zset_{<0}$ $(\b \in \Zset_{>0}),$ $i<|\b|$ by
\eqref{2.27}, \eqref{2.28}. (Note that these cases were excluded 
in \leref{2.4}.)

\bth{2.5}Assuming \eqref{cond} the following relations
\begin{align}
\lb{2.29} 
&(\Labe(n, T)-1) \ph_{+}^{(i)}(n)= \ph_{+}^{(i-1)}(n),
\\
\lb{2.30}
&(\Labe(n, T)-1) \ps_{+}^{(i)}(n)= \ps_{+}^{(i-1)}(n),
\end{align} 
hold for all $i \in \Zset_{\geq 0}$ if $\a \notin \Zset$ and for 
$i=0, \ldots, |\a|-1$ if $\a \in \Zset.$ Similarly one has
\begin{align}
\lb{2.31} 
&(\Labe(n, T)+1) \ph_{-}^{(i)}(n)= \ph_{-}^{(i-1)}(n),
\\
\lb{2.32}
&(\Labe(n, T)+1) \ps_{-}^{(i)}(n)= \ps_{-}^{(i-1)}(n),
\end{align}
for all $i \in \Zset_{\geq 0}$ if $\b \notin \Zset$ and for 
$i=0, \ldots, |\b|-1$
if $\b \in \Zset.$ {\em{(}}We set 
$\ph_\pm^{(-1)}(n) = \ps_\pm^{(-1)}(n) = 0.)$

The kernels of $(\Labe(n, T)-1)^k$ and
$(\Labe(n, T)+1)^l$ are given by
\begin{align}
\lb{2.33} 
&\ker (\Labe(n, T)-1)^k = 
  \spa \{ \ph_{+}^{(i)}(n), \ps_{+}^{(i)}(n) \}_{i=0}^{k-1},
\\
\lb{2.34} 
&\ker (\Labe(n, T)+1)^l =
  \spa \{ \ph_{-}^{(i)}(n), \ps_{-}^{(i)}(n) \}_{i=0}^{l-1},
\end{align}
for $k \leq |\a|$ if $\a \in \Zset,$ for $l \leq |\b|$ if 
$\b \in \Zset ,$ and for all $k, l \geq 0$ if 
$\a, \b \notin \Zset.$ 
\eth
\begin{proof}
In the case $\a \notin \Zset,$ the functions $\ph_{+}(n,z)$ and
$\ps_{+}(n,z)$ are well defined. From the remark in the beginning of this
subsection it follows that \eqref{2.19} and the definitions
of $\ph_{+}^{(j)}(n),$ $\ps_{+}^{(j)}(n)$ imply 
\eqref{2.29}, \eqref{2.30}. The case $\a \in \Zset,$
$i < |\a|$ follows by continuity on $\a.$

The inclusion $\supset$ in \eqref{2.33}, \eqref{2.34} clearly 
follows from \eqref{2.29}--\eqref{2.32}. Because 
$(\Labe(n,T) -1)^k$ is a regular difference operator with 
support $[-k, k],$ to prove \eqref{2.33} it suffices to 
show that the functions $\ph_{+}^{(i)}(n),$ $\ps_{+}^{(i)}(n),$
$i=0, \ldots, k-1$ are linearly independent.

Assume that 
\[
\sum_{i=0}^{k_0} \left( 
         a_i \ph_{+}^{(i)}(n) +
         b_i \ps_{+}^{(i)}(n) \right) =0,
\quad \forall n \in \Zset
\]
for some complex numbers $a_0, \ldots, a_{k_0},$ 
$b_0, \ldots, b_{k_0},$ such that $a_{k_0} \neq 0$ or
$b_{k_0} \neq 0$ $(k_0 \leq k-1).$ Applying $(\Labe(n,T))^{k_0-1}$ to
this
equality and using \eqref{2.29}, \eqref{2.30}, we get
\[
a_{k_0} \ph_{+}^{(0)}(n) + b_{k_0} \ps_{+}^{(0)}(n) = 0, 
\quad \forall n \in \Zset,
\]
i.e.
\begin{equation}
\lb{lindepend}
a_{k_0}\frac{(\e+\a+1)_n}{(\e+1)_n}
= - b_{k_0}\frac{(\e+\b+1)_n}{(\e+\a+\b+1)_n},
\quad \forall n \in \Zset.
\end{equation}
For $n=0$ this gives $a_{k_0}= - b_{k_0} \: ( \neq 0).$
Dividing the two sides of eq. \eqref{lindepend} for two consecutive 
values of $n,$ we get
\[
(\e+\a+n)(\e+\a+\b+n) = (\e+n)(\e+\b+n), \quad \forall n \in \Zset
\]
This gives $\a=0$ which is a contradiction.
Eq. \eqref{2.34} is proved analogously.
\end{proof}

\bre{2.ker}It is clear that
\[
\ker (\Labe -1)^k \cap \ker (\Labe + 1)^l = \emptyset.
\]
Therefore
\[
\ker (\Labe -1)^k (\Labe + 1)^l = 
\ker (\Labe -1)^k \oplus \ker (\Labe + 1)^l 
\]
and \thref{2.5} describes the kernel of the operator
$(\Labe - 1)^k (\Labe + 1)^l$ in the cases specified there.
\ere
\subsection{The sets $\D$ of Darboux transformations from
$\Labe(n,z)$}
Let us fix two nonnegative integers $k$ and $l$ and choose $2(k+l)$
complex numbers
\begin{align*}
& A_i, B_i, \; i=0, \ldots, k-1, 
\\
& C_j, D_j, \; j=0, \ldots, l-1. 
\end{align*}
If $k>0$ $(l>0)$ we will assume $\a \neq -k+1, \ldots,k-1$ 
$(\b \neq -l+1, \ldots, l-1).$ Set
\[
f^{(i)}(n) = 
               \begin{cases}
                  \sum_{r=0}^i(A_r \ph_{+}^{(i-r)}(n)  +
                   B_r \ps_{+}^{(i-r)}(n)),
                   &\text{for $i= 0, \ldots, k-1$} \\
                  \sum_{r=0}^{i-k}(C_r \ph_{-}^{(i-k-r)}(n)  +
                   D_r \ps_{-}^{(i-k-r)}(n)), 
                   &\text{for $\; i= k, \ldots, k+l-1$}
               \end{cases} 
\]
The values of the parameters $A, B, C, D \in \Cset$ for which
\begin{equation}
\lb{cond--det}
\det(n)= \det(f^{(i)}(n+j) )_{i,j=0,-k-l}^{k+l-1, -1}
\neq 0, \; \forall n \in \Zset,
\end{equation}
will be called {\em{admissible.}} For such values we define the
operator 
\begin{equation}
\lb{2.41}
P(n, T) = 
\frac{1}{\det(n)}
\begin{vmatrix}
  f^{(0)}(n-k-l) & \cdots &
  f^{(k+l-1)}(n-k-l) &
  T^{-(k+l)} \cr
  \cdots & \cdots & \cdots &
  \cdots \cr
  f^{(0)}(n) & \cdots &
  f^{(k+l-1)}(n) & 1 
\end{vmatrix}
.
\end{equation}
By expanding \eqref{2.41} along the last column
one sees that the term of $T^{-(k+l)}$ is given by
\[
\frac{\det(n+1)}{\det(n)} \neq 0
\]
hence $P(n, T)$ is a regular difference operator.
As a consequence of properties \eqref{2.29}--\eqref{2.32}  
we obtain
\begin{equation}
\lb{jor1}
(\Labe(n, T)-1) f^{(0)}(n)=0, \quad
(\Labe(n, T)-1) f^{(i)}= f^{(i-1)}(n)
\end{equation}
for $i=1, \ldots, k-1$ and
\begin{equation}
\lb{jor2}
(\Labe(n, T)+1) f^{(k)}(n)=0, \quad
(\Labe(n, T)+1) f^{(j)}= f^{(j-1)}(n)
\end{equation}
for $j=k+1, \ldots, k+l-1.$
Thus $\ker P(n, T)= \spa \{ f^{(i)}(n) \}_{i=0}^{k+l-1}$
is preserved by $\Labe(n,T)$ and according to \prref{darb}
there exists a difference operator $L(n,T)$ with support
$[-1,1]$ such that
\begin{equation}
\lb{defL}
L(n, T) P(n, T) = P(n, T) \Labe(n, T).
\end{equation}
The set of all difference operators $L(n,T)$
for admissible values of the parameters $A, B, C, D$
will be denoted by 
\[
\D.
\]
All operators $L(n, T) \in \D$ are Darboux transformations from 
$\Labe(n, T)$ and $k, l$ refer to the multiplicity of the eigenvalues $1$
and $-1$ of $\Labe(n, T)$ in $\ker P(n,T),$ see eqs. \eqref{jor1} and
\eqref{jor2}. (Recall from part (i) of \prref{darb} that $\Labe(n, T)$
preserves $\ker P(n,T)).$ Every $L(n,T) \in \D$ is a {\em{regular}}
difference operator with eigenfunction 
\begin{equation}
\lb{Psi}
\Ps(n,z) = P(n, T) \pabe(n,z),
\end{equation}
more precisely:
\begin{equation}
\lb{eigenfun}
L(n,T) \Ps(n,z) = z \Ps(n,z).
\end{equation} 

The admissibility condition \eqref{cond--det} holds for almost 
all values of $A, B, C, D \in \Zset.$ 
The complement of the corresponding set in
$\Cset^{2(k+l)}$ consists of the zeros of countably many 
polynomials, obtained from $\det(n)$ 
for fixed $n\in \Zset$ (recall \eqref{cond--det}).
The latter do not vanish identically
due to the linear independence of the set of functions
$\{\ph_\pm^{(i)}(n)\}_{i=0}^{k-1} \cup 
\{ \ps_\pm^{(j)}(n) \}_{j=0}^{l-1}$
(see the proof of \thref{2.5}) and the regularity of $\Labe(n,T).$

There are in fact $k+l$ free parameters in the definition of an
element $L(n, T)\in \D$ since the operator $P(n, T)$ (see \eqref{2.41}) 
only depends on the 
choice of the space $\spa \{ f^{(i)}(n) \}_{i=0}^{k+l-1} (=
\ker P(n,T)),$
and not on the choice of the individual functions
$f^{(i)}(n).$ Using again the linear independence of 
$\{\ph_\pm^{(i)}(n)\}_{i=0}^{k-1} \cup 
\{ \ps_\pm^{(j)}(n) \}_{j=0}^{l-1},$ 
the choice of span is equivalent to a choice of flags
\[
V_0 \sub V_1 \sub \ldots \sub V_{k-1} \; \mbox{and} \; 
W_0 \sub W_1 \sub \ldots \sub W_{l-1}
\]
where $V_i = \spa \{ f^{(r)}(n) \}_{r=0}^{i}$
and $W_j=\spa \{ f^{(r)}(n) \}_{r=k}^{k+j},$ cf. \cite{GHH}.

The relations \eqref{rel1} and \eqref{rel2} for $\Labe(n,T)$ imply
similar relations for the sets $\D:$ 
\begin{align}
\lb{Drel1}
&\Dd_{-\a, -\b, \e+\a+\b}^{(k,l)} = \D, \\
\lb{Drel2}
&\Ad_{(-1)^n} \Dd_{\b, \a, \e}^{(l,k)} =
 \Ad_{(-1)^n} \Dd_{-\b, -\a, \e+\a+\b}^{(l,k)}= - \D.
\end{align}
Here, in addition to \eqref{rel1} and \eqref{rel2}, 
we use that the change of parameters
$\a \ra -\a,$ $\b \ra -\b,$ $\e \ra \e+ \a + \b$ exchanges
$\ph^{(i)}_+(n)$ with $\ps^{(i)}_+(n)$ and  
$\ph^{(i)}_-(n)$ with $\ps^{(i)}_-(n).$ Analogously
the change of parameters 
$\a \ra \b,$ $\b \ra \a,$ $\e \ra \e$
exchanges  
$\ph^{(i)}_+(n)$ with $(-1)^n \ph^{(i)}_-(n)$ and 
$\ps^{(i)}_+(n)$ with $(-1)^n \ps^{(i)}_-(n).$

The Darboux maps between Jacobi functions (operators)
represented by the first identities in \eqref{1.22}, \eqref{1.24}
and \reref{comp} imply the following inclusion relations
\begin{align}
\lb{inclusion}
\Ad_{\frac{2n + 2\e +\a +\b}{n+\e+\a+\b} }
\Dd^{(k-1, l)}_{\a-1, \b;\e} \sub
\D,
\\
\lb{inc}
\Ad_{\frac{2n + 2\e +\a +\b}{n+\e+\a+\b} }
\Dd^{(k, l-1)}_{\a, \b-1;\e} \sub
\D.
\end{align}
The function $(n+\e+\a+\b)/(2n + 2\e +\a +\b)$ is the 
leading coefficient
of the the operators $D^\a_-(n,T)$ and $D^\b_-(n,T),$
see Sect.~2.3. 
Recall that the operator $P(n,T)$ is normalized
to have leading coefficient $1.$
\bre{Jordan}Note that \eqref{jor1}, \eqref{jor2} imply that for the
operator $P(n, T)$ \eqref{2.41} defining an element $L(n,T)$ in $\D$ the
endomorphism $\Labe(n,T)$ on $\ker P(n,T)$ has two Jordan blocks with
eigenvalues $1$ and $-1$ and lengths $k$ and $l,$ respectively.
Insisting on multiple blocks with equal eigenvalues
does not produce larger sets of transformations since the operator 
$\Labe(n,T)$ has a two dimensional kernel. Allowing
$k > |\a|$ or $l > \b$ in the cases $\a \in \Zset$ or $\b \in \Zset$
causes the operators $P(n,T)$ and $L(n,T)$ to have nonrational
coefficients which does lead to bispectrality of $L(n,T)$ as was noted in   
the introduction.
\ere
\hfill \\ 

At the end of this subsection we compute explicitly the coefficients of
the operators $L(n, T)$ in $\D.$ Set
\begin{equation}
\lb{2.47}
L(n, T) = a(n) T + b(n) + c(n) T^{-1}
\end{equation}
for some functions $a(n),$ $b(n),$ and $c(n)$
(the dependence on $A, B, C, D$ will not be shown).
For convenience we denote the coefficients of the
operator $\Labe(n, T)$ by $a_0(n),$ $b_0(n),$ and $c_0(n):$
\begin{equation}
\lb{2.48}
\Labe(n, T) = a_0(n) T + b_0(n) + c_0(n) T^{-1}
\end{equation}
(cf. eq. \eqref{1.17} for their values). Set also
\begin{equation}
\lb{detr}
\det_{-r}(n):=
\det( f^{(i)}(n+j) )_{
\begin{subarray}{l}
i=0, \ldots, k+l-1 \\
j=-k-l, \ldots, -\hat{r}, \ldots, 0
\end{subarray}
} 
\; \mbox{for} \; r=0, \ldots, k+l.
\end{equation}
Note that
\begin{equation}
\lb{determs}
\det_0(n) =  \det(n) \quad \mbox{and} \quad \det_{k+l}(n)=\det_0(n+1)
                                                    =\det(n+1).
\end{equation}
Expanding the determinant \eqref{2.41} defining $P(n,T)$
along the last column gives
\begin{equation}
\lb{2.39}
P(n, T)
= \sum_{r=0}^{k+l} 
(-1)^r \frac{\det_{-r}(n)}{\det(n)}T^{-r}.
\end{equation}

\bpr{2.6}The coefficients $a(n),$ $b(n),$ and $c(n)$ 
of an operator $L(n, T)\in \D$ are expressed in terms of
the coefficients $a_0(n),$ $b_0(n),$ and $c_0(n)$ of
$\Labe(n, T)$ and the functions functions $f^{(i)}(n)$
{\em{(}}see eq. \eqref{detr}{\em{\/)}}
by the following formulas
\begin{align}
\lb{2.52} 
&a(n)=a_0(n),
\\
\lb{2.53}
&b(n)=b_0(n) + a_0(n) \frac{\det_{-1}(n+1)}{\det(n+1)} 
- a_0(n-1) \frac{\det_{-1}(n)}{\det(n)},
\\
\lb{cn}
&c(n)= c_0(n-k -l) \frac{\det(n-1)\det(n+1)}{(\det(n))^2}.
\end{align}
\epr
\begin{proof}
We compare the coefficients of $T$ and $1$ in \eqref{defL}
and use eq. \eqref{2.39} for the operator $P(n, T).$ This
gives the formulas
\begin{align*}
& a(n) = a_0(n),
\\
& b(n) - a(n) \frac{\det_{-1}(n+1)}{\det(n+1)} = 
b_0 (n) - a_0 (n-1) \frac{\det_{-1}(n)}{\det(n)},
\end{align*}
which are equivalent to eqs. \eqref{2.52}, \eqref{2.53}.

Similarly comparing the coefficients of $T^{-k-l-1}$ 
in \eqref{defL} gives
\[
c(n) \frac{\det_{-(k+l)}(n-1)}{\det(n-1)}= 
c_0(n-k-l) \frac{\det_{-(k+l)}(n)}{\det(n)}
\]
which implies \eqref{cn}, taking into account \eqref{determs}. 
\end{proof}
\sectionnew{Bispectral Darboux transformation and an involution}
This section is a preparation for the next one where we show that the
difference operators from $\D$ are bispectral under some natural
conditions on $\a$ and $\b.$ Our proof is based on a result of
\cite{BHY} on Darboux transformations that preserve the bispectral
property. Its application to the situation  under consideration is
nontrivial and requires an intrinsic characterization of a certain space
of difference operators. This is done in terms of an involution of the
algebra of difference operators with rational coefficients.
\subsection{A theorem on bispectral Darboux transformations}
For a fixed choice of the parameters $\a, \b, \e$ we define $\Bn$ 
as the algebra of difference operators $S(n, T)$ with {\em{rational}}
coefficients for which there exists a differential operator $G(z, \pz)$
(also with rational coefficients) such that
\begin{equation}
\lb{3.1}
S(n, T) \pabe(n, z) = G(z, \pz) \pabe(n, z).
\end{equation}
The set of all such operators $G(z, \pz)$ is an algebra which will be
denoted
by $\Bz.$ It is clear that 
\begin{equation}
\lb{3.a}
b \left( S(n, T) \right) := G(z, \pz)
\end{equation}
correctly defines a map
\begin{equation}
\lb{3.2}
b: \Bn \ra \Bz
\end{equation}
which is an antiisomorphism of algebras. 
In this setting eqs. \eqref{1.19}, \eqref{1.20} mean that 
$\lae(n), \Labe(n, T) \in \Bn,$  
$z, \Bab(z, \pz) \in \Bz,$
and
\begin{align}
\lb{3.3} 
&b( \lae(n)) = \Bab(z, \pz),
\\
\lb{3.4}
&b( \Labe(n, T)) = z.
\end{align}
The triple $(\Bn, \Bz, b)$ is an example of a bispectral triple
in the sense of \cite{BHY}.
Denote
\begin{align}
\lb{3.5} 
&\Kn = \Bn \cap \Cset(n),
\\
\lb{3.6} 
&\Kz = \Bz \cap \Cset(z),
\end{align}
where $\Cset(n)$ and $\Cset(z)$ stand for the algebras of rational
functions in the variables $n$ and $z,$ respectively. 
Let
\begin{align}
\lb{3.9} 
&\An = b^{-1} \left( \Kz \right),
\\
\lb{3.10} 
&\Az = b       \left( \Kn \right).
\end{align}
It is obvious that 
\begin{align}
\lb{3.8} 
&\Kz = \Cset [ z ],
\\
\lb{3.12}
&\Az = \Cset[ \Bab(z, \pz) ],
\end{align}
and
\begin{align}
\lb{3.7}
&\Kn \sup \Cset [ \la(n) ],
\\
\lb{3.11} 
&\An \sup \Cset[ \Labe(n, T) ].
\end{align}
Later in \reref{equal} we will show that the inclusions in 
\eqref{3.7} and \eqref{3.11}
can be strengthen to give two equalities.

As was noted in Sect.~3.1, if a difference operator 
$q(\Labe(n, T)) \in \An$
$(q(x) \in \Cset[x])$ is factorized as a product of two 
operators $Q(n, T)$ and $P(n, T)$  
\[
q(\Labe(n, T))=Q(n, T)P(n, T),
\]
then the function 
\[
\Ps(n, z) = P(n, T) \pabe(n, z)
\]
is an eigenfunction of the difference operator $P(n, T) Q(n, T):$
\[
P(n, T) Q(n, T) \Ps(n, z)
= q(z) \Ps(n, z).
\]
We will give a version of Theorem~1.2 from \cite{BHY} 
which provides general sufficient conditions on the operators
$P(n, T)$ and $Q(n, T)$ under which $\Ps(n, z)$ is also an eigenfunction
of a differential operator in the variable $z.$ 
(The original result of \cite{BHY} deals with ``bispectral'' Darboux
transformations in an arbitrary associative algebra but in the form
to be used, needs an additional refinement.)

\bth{3.1}Assume that the operator $q(\Labe(n, T)) \in \An$ is
factorized as
\begin{equation}
\lb{3.13}
q(\Labe(n, T)) = (Q_0(n, T)\nu(n)^{-1}) (\mu(n)^{-1} P_0(n, T))
\end{equation}
for some difference operators $P_0(n, T), Q_0(n, T) \in \Bn$ and
rational functions  $\mu(n), \nu(n) \in \Kn,$ such that the
coefficients of the operators $\mu(n)^{-1} P_0(n, T),$  
$Q_0(n,T) \nu(n)^{-1}$ are correctly defined for $n \in \Zset.$ 
Then the function 
\begin{equation}
\lb{3.14}
\Ps(n, z) = (\mu^{-1}(n) P_0(n, T)) \pabe(n, z)
\end{equation}
satisfies the relations
\begin{align}
\lb{3.15} 
&(\mu(n)^{-1} P_0(n, T)) (Q_0(n, T) \nu(n)^{-1}) \Ps(n, z) =
q(z) \Ps(n, z),
\\
\lb{3.16}
&b(P_0)(z, \pz) b(Q_0)(z, \pz)q(z)^{-1} \Ps(n, z)=
\mu(n) \nu(n) \Ps(n, z),
\end{align}
i.e. it is bispectral.
\eth

Note that in \thref{3.1} we do not assume that the rational functions 
$\mu(n)^{-1}$ and $\nu(n)^{-1}$ are well defined for $n \in \Zset,$ but
only that the ``ratios'' $\mu(n)^{-1} P_0(n, T)$ and
$Q_0(n, T) \nu(n)^{-1}$ are. Because of this a small modification of the
original proof from \cite{BHY} is necessary. 

First of all since 
the algebra $\Bz$ has no zero divisors, eq. \eqref{3.13} implies
(see \cite{BHY}) 
\begin{equation}
\lb{ino}
(b \nu)(z, \pz) \, (b \mu)(z, \pz) = 
(b P_0)(z, \pz) q(z)^{-1} (b Q_0)(z, \pz).
\end{equation}
For all values of $n$ for which $\mu(n)$ does not vanish we have
\[
\Ps(n, z) = \mu(n)^{-1} (b P_0)(z, \pz) \pabe(n, z)
\]
and \eqref{3.16} holds, as a consequence of \eqref{ino}.
The validity of
\eqref{3.16} for all $n \in \Zset$ follows from the definition 
\eqref{3.14} of $\Ps(n, z)$ and the fact that $\pabe(n, z)$ has an
expansion in $z$ around $z=1$ with coefficients that are rational
functions in $n$ (recall \eqref{expansion}).

Returning to the sets $\D$ of Darboux transformations from the operators
$\Labe(n, T),$ we need to find which of the operators
$P(n,T)$ from eq. \eqref{2.41} can be expressed in the form
$\mu(n)^{-1} P_0(n, T)$ with $\mu(n)$ and $P_0(n,T)$ as above. 
According to \thref{3.1} the corresponding
operators $L(n, T) \in \D$ will be bispectral
with bispectral eigenfunction \eqref{Psi} (see also \eqref{3.14}).  
For this we need an invariant description of 
the linear space of difference operators 
\begin{align}
\lb{3.Rn}
\Rn = \spa \{
      \mu(n)^{-1} S(n, T) \mid &
      S(n, T) \in \Bn, \, \mu(n) \in \Kn, \, \mbox{such that} 
\\
\nn   & \mu(n)^{-1} S(n, T) \, \mbox{is well defined for} \, 
      n \in \Zset \}.
\end{align}
This will be obtained in the next subsection. Here we would
like to note that the dual object -- the linear space of differential
operators 
\begin{equation}
\lb{3.Rz}
\Rz = \spa \{
      g(z)^{-1} G(z, \pz) \mid
      G(z, \pz) \in \Bz, \: g(z) \in \Kz
      \}
\end{equation}
is much easier to describe. It is just the space of differential
operators with rational coefficients. This is a consequence of the fact
that the commutator
\[
\left[ \Bab(z, \pz), z \right]
= 2(z^2-1) \pz +((\a-\b) + (\a+\b+2)z)
\]
is a first order differential operator that belongs to $\Bz$ and
$z \in \Kz$ (see eq. \eqref{3.8}). 
Unfortunately for our proof of the fact that the operators from $\D$ are
bispectral we need the space $\Rn,$ and not 
the space $\Rz.$
\subsection{Description of $\Rn$} Denote by $\Del$ the abstract algebra of
difference operators $M(n,T)$ with rational coefficients; that is the
algebra over $\Cset,$ generated by rational functions in $n,$ the shift
operator $T,$ and its inverse $T^{-1},$ subject to the relation
\[
T h(n) = h(n+1) T, \, 
\mbox{for all rational functions} \, h(n).
\] 
Here we do not require that the coefficients of an operator $M(n, T)$
in $\Del$ be well defined for $n \in \Zset.$ More explicitly these
coefficients could have poles at some $n \in \Zset.$
{\em{The subspace of $\Del$ consisting of
operators having this extra regularity property will be denoted by
$\Del^\reg.$}}
We will identify the space of difference operators with
rational coefficients acting on functions $f: \Zset \ra \Cset$ with 
$\Del^\reg.$ In particular, $\Bnt \sub \Rnt \sub \Del^\reg.$ 

Define an involution $I$ in the algebra $\Del$ acting on rational
functions $h(n)$ by
\begin{equation}
\lb{3.17}
(Ih)(n) = h(-(n+ 2 \e + \a + \b +1))  
\end{equation}
and on the shift operator $T$ by
\[
I(T) = T^{-1}.
\]
The involution $I$ is correctly defined since
\[
I(T) \: (Ih)(n) = (Ih)(n+1) \: I(T).
\]
Denote the fixed points of $I$ in $\Del$ by $\Del^I:$
\begin{equation}
\lb{fxI}
\Del^I= \{ M(n, T) \in \Del \mid I(M(n,T))=M(n,T) \}.
\end{equation}

Let
\begin{equation}
\lb{3.vph}
\vph(n)= \frac{(\e+\a+1)_n}{(\e+1)_n} 
\end{equation}
(cf. the definition \eqref{1.16} of $\pabe(n,z)$ for 
$\a \notin \Zset_{<0}).$

\bth{3.2} The space of difference operators $\Rn$ defined in \eqref{3.Rn}
is characterized by
\begin{equation} 
\lb{3.19}
\Rn = \Ad_{\vph(n)} \left( \Del^I \cap \Del^\reg \right),
\end{equation}
i.e. after conjugation by $\vph(n)^{-1}$ all operators from
$\Rnt$ are $I$-invariant. 
\eth
\begin{proof}
Consider first the case $\a \notin \Zset_{<0}.$ Let
\begin{equation}
\lb{3.p}
\pabet(n, z) = \vph(n)^{-1} \pabe(n, z).
\end{equation}
The expression \eqref{1.16} implies
\begin{equation}
\lb{3.20}
\pabet(n, z) = F( -(-n+\e), n + \e +\a + \b+1; \a+1 ; (1-z)/2).
\end{equation}
Let $\Bnt,$ $\Bzt,$ $\Knt,$ $\Kzt,$ $\Rnt,$ and $\Rzt,$ denote the
$\B,$ $\K$ and $\R$ objects associated with the functions $\pabet(n, z)$
(see the beginning of Sect.~4.1 and eqs. \eqref{3.5}, \eqref{3.6},
\eqref{3.Rn}, \eqref{3.Rz} for the appropriate definitions). Obviously
\[
\Rnt = \Ad_{\vph(n)} \Rn, \;
\Bnt = \Ad_{\vph(n)} \Bn,
\]
and $\Knt= \Kn,$ $\Rzt = \Rz,$ $\Kzt = \Kz,$ $\Bzt = \Bz.$ 
In this notation, the statement of the theorem is equivalent to 
\begin{equation}
\lb{star}
\Rnt= \Del^I \cap \Del^\reg.
\end{equation}

To prove that the l.h.s. of \eqref{star} is contained in the r.h.s.,
let us fix an operator $\wt{R}(n, T) \in \Rnt.$ There exists a 
difference operator $\wt{S}(n, T) \in \Bnt$ and a function
$\t{\mu}(n) \in \Knt$ such that
$\wt{R}(n, T) = \t{\mu}(n)^{-1} \wt{S}(n, T).$
We will prove that all operators from $\Bnt$ are $I$-invariant. 
This in particular shows that all functions from $\Knt \sub \Bnt$ are 
$I$-invariant and so are all operators from $\Rnt.$ 

If $\wt{S}(n, T) \in \Bnt,$ then there exists
a differential operator $G(z, \pz)$ for which
\begin{equation}
\lb{3.R1}
\wt{S}(n, T)  \pabet(n, z) =
G(z, \pz)\pabet(n, z).
\end{equation}
The fact that the hypergeometric
function $F(a, b; c;x)$ is symmetric with respect to  $a$ and $b,$
and formula \eqref{3.20} for $\pabet(n, z)$ imply
\begin{equation}
\lb{3.R2}
I \left( \wt{S}(n, T) \right)  \pabet(n, z) =
G(z, \pz)\pabet(n, z).
\end{equation}
Combining \eqref{3.R1} and \eqref{3.R2}, we 
conclude that
\[
\Big( \wt{S}(n, T) - I \left( \wt{S}(n, T) \right) \Big)
\pabet(n, z)=0.
\]
This is only possible if
\[
I \left( \wt{S}(n, T) \right) =
\wt{S}(n, T). 
\]

The harder part of the proof of \eqref{star} is to show that any 
$I$-invariant difference operator from $\Del^\reg$ belongs to $\Rnt.$
It is sufficient to prove that for any $\wt{R}(n, T) \in \Del^I$ there
exists $\wt{S}(n,T) \in \Bnt$ and $\t{\mu}(n) \in \Knt$ such that
\[
\wt{R}(n,T) = \t{\mu}(n)^{-1} \wt{S}(n,T).
\]
First let us write formulas \eqref{1.19} and \eqref{1.20} in terms of
$\pabet(n, z).$ 
Eq. \eqref{1.20} remains unchanged:
\begin{equation}
\lb{3.23}
\lae(n) \pabet(n, z) = \Bab(z, \pz) \pabet(n, z),
\end{equation}
while eq. \eqref{1.19} becomes
\begin{equation}
\lb{3.24}
\Labet(n, T) \pabet(n, z) = z \pabet(n, z)
\end{equation}
with
\begin{align}
\lb{Labet}
&\Labet(n, T)=\vph(n)^{-1} \Labe(n, T) \vph(n)
\\
\nn
&=\frac{2(n+\e+\a+1)(n+\e+\a+\b+1)}{(2n+2\e+\a+\b+1)(2n+2\e+\a+\b+2)}T
\\
\nn
&+\frac{\a^2 - \b^2}{(2n+2\e+\a+\b)(2n+2\e+\a+\b+2)}
\\
\nn
&+\frac{2(n+\e)(n+\e+\b)}{(2n+2\e+\a+\b)(2n+2\e+\a+\b+1)}T^{-1}.
\end{align}
The algebra
$\Del$ has a natural $\Zset_{\geq 0}$ filtration where $\Del_d$
consists of
all operators from $\Del$ with support $[-d,d].$ 
Denote $\Del^I_d = \Del^I \cap \Del_d.$

We will prove that any difference operator $R^I_d(n, T) \in \Del^I_d,$
can be decomposed as a sum
\begin{equation}
\lb{decomp}
R^I_d(n, T)= \t{\mu}(n)^{-1} \wt{S}(n, T) + R^I_{d-1}(n, T)
\end{equation}
where
\begin{align} 
\lb{3.28} 
&\wt{S}(n, T)     \in \Bnt, \, 
\t{\mu}(n) \in \Knt,
\\
\lb{3.29}
&R^I_{d-1}(n, T) \in \Del^I_{d-1}. 
\end{align}
Since $\Del^I_0= \Cset(\la_e(n))$ (any $I$-invariant rational 
function in $n$ is a rational function in $\la_e(n)),$ 
by induction on $d$ eq. \eqref{decomp} implies that  
\[
R^I_d(n, T) \in \Rnt.
\]
A straightforward computation yields
\begin{align}
\lb{comp1}
\ad_{\lae(n)} T^d 
                    &= \left( \lae(n) - \lae(n+d) \right) T^d
\\
                    &= -d(2n +2\e+\a+\b +d+1) T^d,
\nn
\end{align}
and thus
\[
\ad_{\lae(n)}(\ad_{\lae(n)} +1) \Labet =
2(\Id + I)((n+\e+\a+1)(n+\e+\a+\b+1)T).
\]
So
\begin{align}
\lb{comp2}
&\left( \ad_{\lae(n)}(\ad_{\lae(n)} +1) \Labet \right)^d 
\\
\nn
&=2^d(\Id +I)
   \left( \prod_{i=1}^{d}(n+\e+\a+i)(n+\e+\a+\b+i)T^d
   \right) + U_{d-1}. 
\end{align}
for some $U_{d-1} \in \Del^I_{d-1}.$ (Here we use the $I$-invariance 
of $\Labe(n,T).)$ Denote for simplicity
\[
c_d(n) = 2^d \prod_{i=1}^{d} \left( (n+\e+\a+i)(n+\e+\a+\b+i) \right)
\]
and let
\[
R^I_d(n, T) = \sum_{i=-d}^d \frac{a_i(n)}{b_i(n)} T^i, 
\; a_i(n), b_i(n) \in \Cset[n].
\]
Using \eqref{comp1} and \eqref{comp2} we obtain 
\begin{align}
\lb{3.33}
& a_d \left(-\frac{1}{2d}\ad_{\lae(n)} - \frac{\a+\b+d+1}{2}-\e \right)
  b_d \left(\frac{1}{2d}\ad_{\lae(n)} - \frac{\a+\b-d+1}{2}-\e \right)
  \times
\\
\nn 
& \times
  c_d \left(\frac{1}{2d}\ad_{\lae(n)} - \frac{\a+\b-d+1}{2}-\e \right)
  \Big( \ad_{\lae(n)}(\ad_{\lae(n)} +1) \Labet \Big)^d
\\
\nn
& =
  (\Id + I) 
  \left(b_d(n) \: (I b_d)(n) \: c_d(n) \: (I c_d)(n) \:
  \frac{a_d(n)}{b_d(n)}
  T^d \right)+U_{d-1}
\end{align}
for some other $U_{d-1} \in \Del^I_{d-1}.$
There exists a polynomial $q_d(n)$ for which
\[
b_d(n) \: (I b_d)(n) \: c_d(n) \: (I c_d)(n)
= q_d( \lae(n) )
\]
because the polynomial in the l.h.s. is clearly $I$-invariant. 
Denote by $\wt{S}(n,T)$ the difference operator in \eqref{3.33}.
The l.h.s. of \eqref{3.33} implies that $\wt{S}(n,T)$ belongs to $\Bnt$ 
and the r.h.s. implies 
\[
R^I_d(n, T) - (q(\lae(n) )^{-1} \wt{S}(n,T)  \in \Del^I_{d-1}
\]
which completes the proof of \thref{3.2}.
\end{proof}

\bre{equal} Any fuction $\wt{\mu}(n) \in \Knt$ is $I$-invariant and 
therefore is a rational function in $\lae(n).$ 
In fact $\wt{\mu}(n)$ should be a polynomial in $\lae(n).$
Indeed if $\wt{\mu}(n)= p(\lae(n))/q(\lae(n))$ for two polynomials
$p(x), q(x) \in \Cset[x]$ such that $q(x) \not| \: p(x),$ then there
exists a differential operator $G(z, \pz)$ with rational coefficients
such that
\[
G(z, \pz) \pabet(n, z) = \frac{p(\lae(n))}{q(\lae(n))} \pabet(n, z)
\]
which implies
\begin{equation}
\lb{imposs}
p( \Bab(z, \pz) ) = G(z, \pz) q( \Bab(z, \pz) ).
\end{equation}
This is impossible; if $z_0$ is a root of $q(x)$ and $p(x)$
of multiplicities $d_1 > d_2,$ then there exist a holomorphic
function $g(z)$ in a domain $\Om \sub \Cset$ such that
\[
(\Bab(z, \pz) - z_0)^{d_1} g(z)=0 \quad \text{and} \quad 
q(\Bab(z, \pz) - z_0) g(z) \neq 0
\]
which contradicts with \eqref{imposs}.
Since $\Kn = \Knt$ we finally obtain
\begin{align*}
& \Kn = \Cset[\la_e(n)],
\\
&  \Az = \Cset[\Bab(z, \pz)],
\end{align*}
as promised following \eqref{3.7} and \eqref{3.11}.
\ere

\bre{3.r} The second order bispectral differential operators 
of the even case of Duistermaat--Gr\"unbaum's classification 
\cite{DG} are obtained as Darboux transformations from 
the Bessel operators
\[
L_k(x, \px)= \px^2 - \frac{k(k-1)}{x^2}, \, 
k \in \Zset + \frac{1}{2}
\]
in the sense of \eqref{darboux}. 
More precisely for each operator $L(x, \px)$ of this family 
there exists a differential operator with rational
coefficients $P(x, \px)$ such that
\[
L(x, \px) P(x, \px)= P(x, \px) L_k(x, \px).
\]
In addition, the operator $P(x, \px)$ satisfies
\begin{equation}
\lb{Pinvar}
P(x, \px) = P(-x, - \px).
\end{equation}

Let $\Ii$ denote the involution of the algebra of differential operators
with rational coefficients induced by the diffeomorphism $x \go -x$ of
$\Cset$ (i.e. $(\Ii S)(x, \px) = S(-x, - \px)).$ Then \eqref{Pinvar} 
means that $P(x, \px)$ is invariant under $\Ii.$ This gives the 
relation of the approach of this paper via the involution
$I$ and the space $\Rn$ to the construction of \cite{DG}.
\ere
\sectionnew{Bispectrality of $\D$}
In this section we prove our main result: when the parameters 
$\a$ and $\b$ are subject to certain natural integrality conditions,
the difference operators from $\D$ are bispectral. As an example,
for each $L(n,T) \in \Dd^{(2,0)}_{2,0;\e}$ we find a dual differential
operator of order $10.$

The conditions \eqref{cond} on $\a,$ $\b,$ $\e$ are assumed throughout
this section.
\subsection{Proof of the main result}
The conjugation by the function
$\vph(n)$ (see \eqref{3.vph}), used in 
\thref{3.2}, leads us to consider the functions
$\Ph_\pm^{(i)}:= \ph_\pm^{(i)}(n)/\vph(n),$ 
$\Ps_\pm^{(i)}:= \ps_\pm^{(i)}(n)/\vph(n).$ 
Because of eqs. \eqref{2.25}--\eqref{2.28} they are explicitly given
by the formulas
\begin{align}
\lb{3.37} 
& \Ph_+^{(i)}(n)=
  \frac{(-(n+\e))_i (n+\e+\a+\b+1)_i}{(\a+1)_i (-2)^i},
\\
\lb{3.38} 
& \Ps_+^{(i)}(n)=
  \frac{(\e+\b+1)_n (\e+1)_n}{(\e+\a+1)_n (\e+\a+\b+1)_n} 
  \frac{(-(n+\e+\a+\b))_i (n+\e+1)_i}{(-\a+1)_i (-2)^i}, 
\\
\lb{3.39}  
& \Ph_-^{(i)}(n)=
  \frac{(\e+\b+1)_n}{(-1)^n(\e+\a+1)_n} 
  \frac{(-(n+\e))_i (n+\e+\a+\b+1)_i}{(\b+1)_i 2^i},
\\
\lb{3.40}  
& \Ps_-^{(i)}(n)=
  \frac{(\e+1)_n}{(-1)^n (\e+\a+\b+1)_n} 
  \frac{(-(n+\e+\a+\b))_i (n+\e+1)_i}{(-\b+1)_i 2^i}. 
\end{align}

\ble{3.3}If $\a \in \Zset,$ then for $i \leq |\a|-1,$
$\Ph_+^{(i)}(n)$ and $\Ps_+^{(i)}(n)$ are $I$-invariant
rational functions of $n.$

If $\a \in \Zset$ and $\b \in \Zset$ then
for $i \leq |\a| -1,$ $j \leq |\b| -1,$
$\Ph_+^{(i)}(n),$ $\Ps_+^{(i)}(n),$ 
$(-1)^n \Ph_-^{(j)}(n),$ and $(-1)^n \Ps_-^{(j)}(n)$
are rational functions of $n,$  
$\Ph_+^{(i)}(n),$ $\Ps_+^{(i)}(n)$ are $I$-invariant, and
\begin{align}
\lb{3.41}
& I \left( (-1)^n \Ph_-^{(j)}(n) \right)=
    (-1)^{\a+\b} \left( (-1)^n \Ph_-^{(j)}(n) \right),
\\
\lb{3.42}
& I \left( (-1)^n \Ps_-^{(j)}(n) \right) =
    (-1)^{\a+\b} \left( (-1)^n \Ps_-^{(j)}(n) \right).
\end{align}
\ele
\begin{proof}
First note that
\begin{align*}
&(-(n+\e))_i (n+\e+\a+\b+1)_i
\\
&=\prod_{r=0}^{i-1}(-(n+\e)+r)(n+\e+\a+\b+1+r)
\\
&= (-1)^k \prod_{r=0}^{i-1}(\la(n)-r(\a+\b+1+r))
\end{align*}
and similarly
\[
(-(n+\e+\a+\b))_i (n+\e+1)_i
=(-1)^i \prod_{r=0}^{i-1}(\la(n)-(\a+\b-r)(r+1))
\]
are $I$-invariant polynomials in $n.$

To prove the first statement of the lemma we
use a similar computation. Restricting to the case 
$\a \in \Zset_{>0}:$
\begin{align*}
&\frac{(\e+\b+1)_n (\e+1)_n}{(\e+\a+1)_n (\e+\a+\b+1)_n} 
\\
&=\frac{(\e+1)_\a (\e+\b+1)_\a}
       {(n+\e+1)_\a (n+\e+\b+1)_\a} 
\\
&=\frac{(\e+1)_\a (\e+\b+1)_\a}
       {\prod_{r=1}^\a (n+\e+r)(n+\e+\b+\a+1-r)} 
\\
&=\frac{(\e+1)_\a (\e+\b+1)_\a}
       {\prod_{r=1}^\a (\la(n)+r(\a+\b+1-r))}.
\end{align*}

The proof of the second statement is analogous. Assuming $\a,$
$\b \in \Zset_{>0}$ and $\b \geq \a$ we obtain
\begin{align*}
& \frac{(\e+\b+1)_n}{(\e+\a+1)_n}=
  \frac{(\e+\a+n+1)_{\b-\a}}{(\e+\a+1)_{\b-\a}}
\\
&=\frac{q(n) \prod_{r=1}^{\left[\frac{\b-\a}2\right]}
          (n+\e+\a+r)(n+\e+\b+1-r)}
       {(\e+\a+1)_{\b-\a}}
\\
&=\frac{q(n) \prod_{r=1}^{\left[\frac{\b-\a}2\right]}
              (\la(n)+(\a+r)(\b+1-r))}
              {(\e+\a+1)_{\b-\a}}
\end{align*}
with
\[
q(n)=  
    \begin{cases}
    1,  &\text{if $\b + \a$ is even} \\
    n + \e +(\a+\b+1)/2, &\mbox{if $\b + \a$ is odd}
    \end{cases}
\]
(Since $\a \in \Zset,$ the first condition is equivalent to $2 | (\b-\a)$
and the second one to $2 \not| (\b - \a).)$ 
To finish the proof of \eqref{3.41} we just observe that
\begin{equation}
\lb{3.I}
I(n + \e +(\a+\b+1)/2) = -(n + \e +(\a+\b+1)/2).
\end{equation}
The remaining cases for $\a, \b \in \Zset$ are treated analogously. 

The identity \eqref{3.42} follows from the analogous formula
\[
\frac{(\e+1)_n}{(\e+\a+\b+1)_n}
=\frac{(\e+1)_{\a+\b}}
       { q(n) \prod_{r=1}^{\left[\frac{\b+\a}{2}\right]}
          (\la(n)+r(\a+\b+1-r))                       }
\]
and eq. \eqref{3.I}.

Throughout this proof, for a real number $x$ by $[x]$ we denote its
integer part.
\end{proof}

\bth{3.4} Assuming \eqref{cond},
the following sets consist of bispectral difference operators:

1.) $\Dko$ if $\a \in \Zset$ and $k \leq |\a|,$

2.) $\Dol$ if $\b \in \Zset$ and $l \leq |\b|,$

3.) $\D$ if $\a, \; \b \in \Zset$ and $k \leq |\a|,$
                                        $l \leq |\b|.$
\eth

When the conditions \eqref{cond} are not met but the operator $\Labe(n,T)$
is still well defined the arguments below can be adapted properly. We 
do not pursue that here.

\begin{proof} Because of the relation \eqref{Drel2} the second case 
follows from the first one.  

Let us restrict to instances 1 and 3 of the theorem above.
In each of them we can assume that
$k+l$ is even using \eqref{inclusion}.
Fix an operator $L(n, z) \in \D,$ determined by
a choice of the functions $\{f^{(i)}(n)\}_{i=0}^{k+l-1},$
i.e. a choice of admissible values of the complex parameters
$A, B,C,D$ (see Sect.~3.3). It has the eigenfunction
$\Ps(n,z)$ defined in \eqref{Psi}
\[
L(n,T) \Ps(n,z) = z \Ps(n,z),
\]
cf. \eqref{eigenfun}. We need to show that there exists
a differential operator $B(z, \pz)$ having $\Ps(n,z)$ as an
eigenfunction, that is
\[
B(z, \pz) \Ps(n,z) = \th(n) \Ps(n,z)
\]
for some function $\th(n).$

Define the functions
\[
F^{(i)}(n)= f^{(i)}(n)/\vph(n), \quad i=0, \ldots, k+l-1.
\]
Let us put $s:=(k+l)/2$ and consider the operator
\begin{equation}
\lb{tP}
\widetilde{P}(n, T)= (-1)^{nl}
\begin{vmatrix}
     F^{(0)}(n-s) & \ldots & F^{(k+l-1)}(n-s) & T^{-s} \cr
     \ldots & \ldots & \ldots & \ldots \cr
     F^{(0)}(n+s) & \ldots & F^{(k+l-1)}(n+s) & T^{s} \cr
\end{vmatrix}
.
\end{equation}
It is a regular difference operator with
kernel given by $\spa \{F^{(i)}(n)\}_{i=0}^{k+l-1}.$ 
Hence it is related to the operator $P(n,T)$ 
(recall eq. \eqref{2.41}) by
\begin{equation}
\lb{3.51}
P(n, T) =
d(n)^{-1} \vph(n)^{-1} T^{-s}
\widetilde{P}(n, T) \vph(n)
\end{equation}
where
\[
d(n) = (-1)^{nl}\det(F^{(i)}(n+j))_{i,j=0, -k-l}^{k+l-1, -1}
\]
is the leading coefficient of $T^{-s}\widetilde{P}(n, T).$
\leref{3.3} implies that $F^{(0)}(n),$ $\ldots,$ $F^{(k-1)}(n),$
and $(-1)^n F^{(k)}(n),$ $\ldots,$ $(-1)^n F^{(k+l-1)}(n)$ 
are rational functions in $n.$ This implies that 
for $i=k, \ldots, k+l-1$ and for all $j \in \Zset,$
$(-1)^n F^{(i)}(n+j)$
are also rational functions in $n$
and thus $\widetilde{P}(n, T)$ has rational coefficients. 
In addition \leref{3.3} gives
\[
I\left(F^{(i)}(n+j)\right) = F^{(i)}(n-j),
\, i= 0, \ldots, k-1, \, j \in \Zset
\]
and
\[
I\left((-1)^n F^{(i)}(n+j)\right) = (-1)^{(\a+\b)/2} 
\left( (-1)^n F^{(i)}(n-j)\right), 
\, i= k, \ldots, k+l-1, \, j \in \Zset.
\]
Taking into account that $I(T)= T^{-1}$ we obtain
\begin{equation}
\lb{3.54}
I \left( \widetilde{P}(n,T) \right) =            
(-1)^{s} (-1)^{(\a+\b)l}
\widetilde{P}(n,T)
\end{equation}
where the factor $(-1)^{s}$ comes from exchanging the pairs of rows
$(1, k+l+1),$ $\ldots,$ $(s, s +2).$
Set
\[
q(n)          = 
              \begin{cases}
                     1, &\text{if $s + (\a+\b)l$ is even} \\
                     (n+\e + (\a+\b+1)/2), 
                       &\text{if $s + (\a+\b)l$ is odd}
              \end{cases}
\]
and consider the operator 
\begin{equation}
\lb{3.55}
\ol{P}(n, T) = q(n) \widetilde{P}(n, T).
\end{equation}
Because of the conditions \eqref{cond}, $q(n)$ does not vanish
for $n \in \Zset.$ Taking into account \eqref{3.51} one sees
that $\ol{P}(n, T)$ is related to $P(n, T)$ by
\begin{equation}
\lb{3.Pbar}
P(n, T) =
\frac{d(n) \vph(n-s)}{q(n-s) \vph(n)}
T^{-s} \vph(n)^{-1} \ol{P}(n,T) \vph(n).
\end{equation}
Since $\ol{P}(n, T)$ is a regular difference operator
and $\vph(n)$ does not vanish for $n \in \Zset$ 
(recall \eqref{cond}),
there exists a difference operator with rational
coefficients $\ol{Q}(n, T)$ such that
\[
(\Labe(n, T)-1)^k (\Labe(n, T)+1)^l =
\left( \vph(n)^{-1} \ol{Q}(n, T) \vph(n) \right)
\left( \vph(n)^{-1} \ol{P}(n, T) \vph(n) \right)
\]

  {} From eqs. \eqref{3.54} and \eqref{3.55} it follows that
$\ol{P}(n,T)$ is $I$-invariant. Finally combining this with
the $I$-invariance of 
$\vph(n)^{-1} \Labe(n, T) \vph(n)=\Labet(n,T)$
(see \eqref{Labet}) implies the 
$I$-invariance of the operator $\ol{Q}(n, T).$ 
\thref{3.2} now gives 
\[
\vph(n)^{-1} \ol{P}(n, T) \vph(n), \, 
\vph(n)^{-1} \ol{Q}(n, T) \vph(n) \in \Rn.
\] 
Applying \thref{3.1} we obtain that the function 
\begin{equation}
\lb{olPs}
\ol{\Ps}(n, z) = \vph(n) \ol{P}(n, T) \vph(n)^{-1} \pabe(n, z)
\end{equation}
is an eigenfunction of a differential operator $B(z, \pz)$
\begin{equation}
\lb{eqolPs}
B(z, \pz) \ol{\Ps}(n, z) = h(\lae(n)) \ol{\Ps}(n, z),
\end{equation}
for some polynomial $h(x).$ Because of \eqref{3.Pbar}
our original function $\Ps(n, z) \in \D$ 
is related to $\ol{\Ps}(n, z)$ by
\begin{equation}
\lb{Psirel}
\Ps(n,z)= P(n, T) \pabe(n, z) =
\frac{d(n) \vph(n-s)}{q(n-s) \vph(n)}
T^{-(k+l)/2} \ol{\Ps}(n, z).
\end{equation}
Eq. \eqref{eqolPs} implies that
$\Ps(n,z)$ is an eigenfunction of the same operator
$B(z, \pz)$ with eigenvalue $T^{-(k+l)/2} h(\lae(n)):$
\begin{equation}
\lb{BPs}
B(z, \pz) \Ps(n, z) = h(\lae(n-(k+l)/2)) \Ps(n, z).
\end{equation}
\end{proof}
\subsection{An example: the set $\Dd_{2, 0, \e}^{(2, 0)}$}
In this final subsection we consider in detail the case
$\a=2,$ $\b=0,$ $k=2,$ $l=0$ and use this example for two
different purposes. First we give the reader a guided tour
through the results in this paper:
we start with the function $p_e^{2,0}(n,z)$ from \eqref{1.16}, 
give the ingredients needed to build the difference operator 
$P(n,T)$ \eqref{2.41} and the corresponding eigenfunction 
$\Psi(n,z)$ \eqref{Psi}, and end with a description of 
the strategy used in the construction of a 
differential operator in the variable $z$ giving a bispectral
situation. The algebra of possible differential operators
in $z$ contains some whose order is 
lower than the one resulting from 
this construction. We close this subsection with an 
explicit expression for the (essentially unique) bispectral 
operator of minimal order and material related to this operator.

The functions $\ph^{(i)}_+(n)$ and $\ps^{(i)}_+(n)$ 
$(i=0,1)$ from \eqref{2.25}--\eqref{2.26} are given by
\begin{align}
\lb{ph}
&\ph^{(0)}_+(n)=\frac{(n+\e+1)_2}{\ka}, \quad
\ph^{(1)}_+(n)=\frac{(n+\e)_4}{6\ka}, 
\\
\lb{ps}
&\ps^{(0)}_+(n)=\frac{\ka}{(n+\e+1)_2}, \quad
\ps^{(1)}_+(n)=-\frac{\ka}{2},
\end{align}
where
\begin{equation}
\lb{ka}
\ka=(\e+1)(\e+2).
\end{equation}
The conditions \eqref{cond} reduce to
$\e \notin \Zset.$

An element $L(n,T) \in \Dd_{2, 0, \e}^{(2, 0)}$ is determined by a choice
of the functions
\begin{align*}
&f^{(0)}(n)= A_0 \ph^{(0)}_+(n) + B_0 \ps^{(0)}_+(n),
\\
&f^{(1)}(n)= A_1 \ph^{(0)}_+(n) + B_1 \ps^{(0)}_+(n) +
             A_0 \ph^{(1)}_+(n) + B_0 \ps^{(1)}_+(n),
\end{align*}
cf. Sect.~3.3.
We will restrict to the generic case when $A_0 \neq 0.$ 
In this case we can assume that $A_0=1$ and $A_1=0$ by 
dividing $f^{(0)}(n)$ by $A_0$ and then subtracting from $f^{(1)}(n)$
the term $A_1 f^{(0)}(n).$ Recall that $L(n, T)$ depends only on
$\spa \{f^{(0)}(n), f^{(1)}(n) \}.$ Once this space has been
specified by the choice of $B_0,$ $B_1$ we can build the difference
operator $P(n,T)$ as in \eqref{2.41} and we get the eigenfunction 
$\Psi(n,z)$ of $L(n, T)$ from \eqref{Psi}.

The theory developed in Sections 4 and 5 makes it convenient to introduce
the difference operators $\wt{P}(n,T),$ see \eqref{tP}, and $\ol{P}(n,T),$
see \eqref{3.55}, related to $P(n,T)$ by \eqref{3.51} and \eqref{3.Pbar}.

The main point in the proof of \thref{3.4} is that the operator
$\ol{P}(n, T)$ defined in \eqref{3.55} (see also \eqref{tP}) 
is $I$-invariant  and thus $\vph(n) \ol{P}(n, T) \vph(n)^{-1} \in 
{{\mathcal R}_{2, 0; \e}}.$ This implies that 
the function
\[
\ol{\Ps}(n, z) = \vph(n) \ol{P}(n, T) \vph(n)^{-1} \pabe(n, z)
\]
(see \eqref{olPs})
can be expressed as
\begin{equation}
\lb{dual}
\ol{\Ps}(n, z) = \mu(n)^{-1} G(z, \pz) \pabe(n, z)
\end{equation}
for some differential operator with rational coefficients $G(z, \pz)$
and some polynomial $\mu(n)$ 
(recall the definition \eqref{3.Rn} of ${{\mathcal R}_{2, 0; \e}}).$
Now any operator $B(z, \pz)$ that is a Darboux transformation from
$h(B_{2,0}(z, \pz))$ for some $h(x) \in \Cset[x]$ via
the operator $G(z, \pz),$ i.e.
\begin{equation}
\lb{zDarboux}
B(z, \pz) G(z, \pz) = G(z, \pz) h(B_{2,0}(z, \pz))
\end{equation}
will satisfy
\[
B(z, \pz) \ol{\Ps}(n, z) = h( \lae(n)) \ol{\Ps}(n, z)
\]
(a differential analog of \eqref{leig}). The function $\Ps(n,z)$ 
is related to $\ol{\Ps}(n,z)$ by \eqref{Psirel} and is
also an eigenfunction of $B(z, \pz)$ but with eigenvalue
$h(\la(n-1))$
\[
B(z, \pz) \Ps(n, z) = h( \lae(n-1)) \Ps(n, z),
\]
see \eqref{BPs}. Combined with \eqref{eigenfun}
\[
L(n, T) \Ps(n, z)= z \Ps(n, z)
\]
this gives the desired bispectral pair $(L(n, T), B(z, \pz)).$

The $I$-invariance of the operator $\ol{P}(n,T)$ in this special case 
can be observed directly. Because of \eqref{ph}, \eqref{ps} the functions
$F^{(i)}(n)= f^{(i)}(n)/\vph(n),$ $i=0,1,$ see \eqref{3.vph}, 
are given in terms of
\[
\lae(n)=(n+\e)(n+\e+3)
\]
by 
\begin{align*}
&F^{(0)}(n)= 1 + \frac{B_0 \lae(n)}{6 \ka},
\\
&F^{(1)}(n)= \frac{B_1 \lae(n)}{6 \ka} + 
   \frac{\ka}{(\lae(n)+1)} 
   \left( \frac{\ka}{(\lae(n)+1)}- \frac{B_0}{2} \right).
\end{align*}
The operator $\ol{P}(n,T)$ is given by
\[
\ol{P}(n,T)= 
(n +\e + 3/2)
\begin{vmatrix}
     F^{(0)}(n-1) & F^{(1)}(n-1) & T^{-1} \cr
     F^{(0)}(n)   & F^{(1)}(n)   & 1 \cr
     F^{(0)}(n+1) & F^{(1)}(n+1) & T \cr
\end{vmatrix}
\]
and it is $I$-invariant because of the $I$-invariance of $\lae(n)$ 
and the skew invariance of the factor in front compensating
the effect of the exchange of first and third row. 
An operator $G(z, \pz)$ satisfying \eqref{dual} is generated from the
proof of \thref{3.2}. It is of high order and the one of minimal order
$10$ has the following form
\begin{align*}
G(z, \pz) &= (z-1)^6 (z+1)^5 \pz^{10} +
(z-1)^5 (z+1)^4 (57z +7) \pz^9+
\\
&+4 (z-1)^4 (z+1)^3 (311 z^2 + 68 z - 43) \pz^8 + 
(3 B_0 \ka^2 (z-1)^2 (z+1)^2+
\\
&+2 (18793 z^4 + 5796 z^3 
    - 15734 z^2 - 3636 z + 1501)) \pz^7+ \ldots 
\end{align*}
\thref{3.1} guarantees that \eqref{zDarboux} is satisfied for 
some polynomial $h(x).$ It also generates such a polynomial but it is
again of high order. The one of minimal order $5$ is given by
\begin{align*}
h(x-2)&= x^5 -5 x^4 +(10B_0 \ka^2 +8) x^3 -
\\
&-(30B_1 \ka^2 + 20B_0 \ka^2 +4) x^2 -15 B_0^2 \ka^4 x.
\end{align*}
Given $G(z, \pz)$ eq. \eqref{zDarboux} determines the dual bispectral
operator $B(z, \pz)$ of $L(n, T)$ of minimal order uniquely. 
It is given by 
\begin{align*}
&B(z, \pz) = (z-1)^5 (z+1)^5 \pz^{10} 
+ 50 (z-1)^4 z (z+1)^4 \pz^9
+5 (z-1)^3 (z+1)^3 \times
\\
& \quad
\times (11 z -5) 
(17z+7) \pz^8 + 160 (z-1)^2 (z+1)^2 (52 z^3 - 7 z^2 -28 z +1) \pz^7
+
\\
& \quad
+(30 B_0^2 \ka^4 z 
   +120 B_1 \ka^2 (z-1) + 120 B_0 \ka^2)\pz^6 
+ (180 B_0 \ka^2 (z-1)^2z(z+1)^2 +
\\
& \quad
+240 (z-1)^2(337 z^3 +504 z^2 + 141 z -30)) \pz^5 +
(-30B_1 \ka^2 (z-1)^2 (z+1)^2 +
\\
& \quad
+ 120 B_0 \ka^2 (z-1)(z+1)(8 z^2-z-3) 
+120(z-1)^2(641 z^2 + 758 z + 161) ) \pz^4 +
\\
& \quad
+(-240 B_1 \ka^2 (z-1) z (z+1)+ 
240 B_0 \ka^2( 7 z^3 -3 z^2 -7 z +1) + 960(z-1)^2 \times
\\
&  \quad
\times (26z +19) ) \pz^3 
+ ( -60 B_1 \ka^2 (z-1) (7z+5) + 120 B_0 \ka^2 (2 z+1) 
(3z-5)+ 
\\
& \quad
+ 1440 (z-1)^2)\pz^2 
- (30 B_0^2 \ka^4 z + 120 B_1 \ka^2 (z-1) + 120 B_0 \ka^2) \pz.
\end{align*}

In the cases $k=1,$ $l=0,1$ and $\e=0$ the dual bispectral operator of
minimal order was determined in \cite{KK, Zh}.  
    

\begin{thebibliography}{AFMO}
    \bibitem{BHYc} B. Bakalov, E. Horozov, and M. Yakimov, {\em{Bispectral
algebras of commuting ordinary differential
operators,}} Comm. Math. Phys. {\bf{190(2)}} (1997), 331--373. 
    \bibitem{BHYd} B. Bakalov, E. Horozov, and M. Yakimov, {\em{Highest
weight modules over the $W_{1+\infty}$ algebra and the bispectral
problem,}} Duke Math. J. {\bf{93}(1)}  (1998), 41--72. 
    \bibitem{BHY} B. Bakalov, E. Horozov, and M. Yakimov, {\em{General
methods for constructing bispectral operators,}} Phys. Lett. A 
{\bf{222(1-2)}} (1996), 59--66.  
    \bibitem{B} Yu. Berest, {\em{Huygens' principle and the bispectral
problem,}} In: The bispectral problem (Montreal, PQ, 1997), 11--30, CRM
Proc. Lecture Notes, 14, Amer. Math. Soc., Providence, RI, 1998.
    \bibitem{BW} Yu. Berest and G. Wilson, {\em{Classification of rings
of differential operators on affine curves,}} IMRN 1999 {\bf{(2)}},
105--109.
    \bibitem{DG} J. J. Duistermaat and F. A. Gr\"unbaum,
{\em{Differential equations in the spectral
parameter,}} Comm. Math. Phys. {\bf{103(2)}} (1986), 177--240. 
    \bibitem{EV} P. Etingof and A. Varchenko, {\em{Traces of intertwiners
for quantum groups and difference equations, I,}} preprint
math.QA/9907181.
    \bibitem{FMTV} G. Felder, Y. Markov, V. Tarasov, and A. Varchenko,
{\em{Differential Equations Compatible with KZ Equations,}} preprint
math.QA/0001184.
    \bibitem{GH} F. A. Gr\"unbaum and L. Haine, {\em{A theorem of Bochner,
revisited,}} Algebraic aspects of integrable systems, 143--172,
Progr. Nonlinear Differential Equations Appl., 26, Birkhäuser Boston,
Boston, MA, 1997. 
    \bibitem{GH2} F. A. Gr\"unbaum and L. Haine, {\em{Associated
polynomials, spectral matrices and the bispectral problem,}} 
Methods and Applications of Analysis {\bf{6(2)}} (1999), 209--224. 
    \bibitem{GHH} F. A. Gr\"unbaum, L. Haine, and E. Horozov, {\em{Some
functions that generalize the Krall-Laguerre polynomials,}} 
J. Comput. Appl. Math. {\bf{106(2)}} (1999), 271--297.
    \bibitem{H} L. Haine, {\em{Beyond the classical orthogonal
polynomials,}} The bispectral problem (Montreal, PQ, 1997), 47--65, CRM
Proc. Lecture Notes, 14, Amer. Math. Soc., Providence, RI, 1998.  
    \bibitem{H2} L. Haine, {\em{The Bochner--Krall problem: some new
perspectives,}} to appear in Proceedings of the NATO workshop 
on Special functions, Tempe, Arizona, 2000.
    \bibitem{HI} L. Haine and P. Iliev, {\em{Commutative rings of
difference operators and an adelic flag manifold,}} IMRN 2000 {\bf{(6)}}, 
281--323.
    \bibitem{HI2} L. Haine and P. Iliev, {\em{A rational analog 
of the  Krall polynomials,}} J. of Phys. A: Math. Gen., to appear.
    \bibitem{HM} E. Horozov and T. Milanov, {\em{Fuchsian bispectral
operators,}} preprint 2000.
    \bibitem{KR} A. Kasman and M. Rothstein, {\em{Bispectral Darboux
transformations: the generalized Airy case,}} Phys. D
{\bf{102(3-4)}} (1997), 159--176. 
    \bibitem{KK} J. Koekoek and R. Koekoek, {\em{Differential equations
for generalized Jacobi polynomials,}} J. Comput. Appl. Math., to appear.
    \bibitem{HLK} H.~L.~Krall, {\em{Certain differential equations
for Tchebicheff polynomials,}} Duke Math.~J. {\bf{4}} (1938),
705--718.
    \bibitem{MOS} W. Magnus, F. Oberhettinger, and R. Soni,
{\em{Formulas and theorems for the special functions of mathematical
physics,}} Springer Verlag, New York, 1966. 
    \bibitem{Mum} D. Mumford, {\em{An algebro--geometric construction of
commuting operators and of solutions to the Toda lattice equation,
Korteweg deVries equation and related nonlinear equation,}} Proceedings of
the International Symposium on Algebraic Geometry (Kyoto Univ., Kyoto,
1977), pp. 115--153, Kinokuniya Book Store, Tokyo, 1978.
    \bibitem{vMM} P. van Moerbeke and D. Mumford, {\em{The
spectrum of difference operators and algebraic curves,}} Acta Math. 
{\bf{143(1-2)}} (1979), 93--154. 
    \bibitem{Wil} G. Wilson, {\em{Bispectral commutative ordinary
differential operators,}} J. Reine Angew. Math. {\bf{442}} (1993),
177--204. 
    \bibitem{Wcm} G. Wilson, {\em{Collisions of Calogero-Moser particles
and an adelic Grassmannian,}} With an appendix by I. G. Macdonald. 
Invent. Math. {\bf{133(1)}} (1998), 1--41. 
    \bibitem{W} P. Wright, {\em{Darboux transformations, algebraic
subvarieties of Grassmann manifolds, commuting flows and bispectrality,}}
Ph.D. Thesis, Univ. California, Berkeley, 1987.
    \bibitem{Zh} A. Zhedanov, {\em{A method of constructing Krall's 
polynomials,}} J. Comput. Appl. Math. {\bf{107}} (1999), 1--20. 
    \end{thebibliography}
\end{document}